%% file: ms.tex
\documentclass[11pt]{article}
\pdfoutput=1

\usepackage[inner=2.25cm,outer=2.25cm,top=3cm,bottom=3cm,bindingoffset=0cm,verbose]{geometry}  %Defines margins etc.
\usepackage[utf8]{inputenc} % allow utf-8 input
\usepackage[T1]{fontenc}    % use 8-bit T1 fonts
\usepackage[affil-it]{authblk}
\usepackage{kpfonts}        % changes the font
\usepackage{booktabs}       % professional-quality tables
\usepackage{amsfonts,amssymb,amsmath,amsthm,amstext,amssymb,amsopn,mathtools,nicefrac,xfrac,bm}
\usepackage{microtype}      % microtypography
\usepackage{lipsum}
\usepackage{graphics,graphicx}
\usepackage[font=small,labelfont={small,bf}]{caption}
\usepackage{subcaption}
\usepackage{epstopdf}
\usepackage{verbatim}
\usepackage{appendix}
\usepackage{array}

\usepackage[usenames,dvipsnames]{xcolor}
\usepackage{hyperref}                            % hyperlinks
\usepackage{url}            % simple URL typesetting
\hypersetup{colorlinks=true,linkcolor=blue,filecolor=magenta,urlcolor=Blue,citecolor=OliveGreen}
\usepackage[noabbrev,capitalise,nosort,nameinlink]{cleveref}       % uses cref instead of ref

\usepackage{fancyhdr}
\usepackage{algorithm,algorithmic}

\providecommand{\keywords}[1]
{
  \small	
  \textbf{Keywords.} #1
}

\title{Stochastic parareal: an application of probabilistic methods to time-parallelisation}

\author[1]{Kamran Pentland\thanks{Corresponding author: \href{mailto:kamran.pentland@warwick.ac.uk}{kamran.pentland@warwick.ac.uk}}}
\author[2]{Massimiliano Tamborrino}
\author[3]{D. Samaddar}
\author[3]{L. C. Appel}
\affil[1]{Mathematics Institute, University of Warwick, Coventry, CV4 7AL, United Kingdom}
\affil[2]{Department of Statistics, University of Warwick, Coventry, CV4 7AL, United Kingdom}
\affil[3]{Culham Centre for Fusion Energy, Culham Science Centre, Abingdon, Oxfordshire, OX14 3DB, United Kingdom}

\date{\today}

\begin{document}
\maketitle

\begin{abstract}
Parareal is a well-studied algorithm for numerically integrating systems of time-dependent differential equations by parallelising the temporal domain. Given approximate initial values at each temporal sub-interval, the algorithm locates a solution in a fixed number of iterations using a predictor-corrector, stopping once a tolerance is met. This iterative process combines solutions located by inexpensive (coarse resolution) and expensive (fine resolution) numerical integrators. In this paper, we introduce a \textit{stochastic parareal} algorithm aimed at accelerating the convergence of the deterministic parareal algorithm. Instead of providing the predictor-corrector with a deterministically located set of initial values, the stochastic algorithm samples initial values from dynamically varying probability distributions in each temporal sub-interval. All samples are then propagated in parallel using the expensive integrator. The set of sampled initial values yielding the most continuous (smoothest) trajectory across consecutive sub-intervals are fed into the predictor-corrector, converging in fewer iterations than the deterministic algorithm with a given probability. The performance of the stochastic algorithm, implemented using various probability distributions, is illustrated on low-dimensional systems of ordinary differential equations (ODEs). We provide numerical evidence that when the number of sampled initial values is large enough, stochastic parareal converges almost certainly in fewer iterations than the deterministic algorithm, maintaining solution accuracy. Given its stochastic nature, we also highlight that multiple simulations of stochastic parareal return a distribution of solutions that can represent a measure of uncertainty over the ODE solution.
\end{abstract}

% keywords can be removed

%%%%%%%%%%%%%%%%%%%%%%%%%%%%%%%%%%%%%%%%%%%%%%%%%%%%%%%%%%%%%%%%%
% Paper content

\input{sec1}
\input{sec2}
\input{sec3}
\input{sec4}

\input{sec5}

%%%%%%%%%%%%%%%%%%%%%%%%%%%%%%%%%%%%%%%%%%%%%%%%%%%%%%%%%
%Acknowledgements
\section*{Acknowledgements}
This work has partly been carried out within the framework of the EUROfusion Consortium and has received funding from the Euratom research and training programme 2014-2018 and 2019-2020 under grant agreement No 633053. The views and opinions expressed herein do not necessarily reflect those of the European Commission. KP is funded by the Engineering and Physical Sciences Research Council through the MathSys II CDT (grant EP/S022244/1) as well as the Culham Centre for Fusion Energy. The authors would also like to acknowledge the University of Warwick Scientific Computing Research Technology Platform for assistance in the research described in this paper.

%%%%%%%%%%%%%%%%%%%%%%%%%%%%%%%%%%%%%%%%%%%%%%%%%%%%%%%%%%%%%%%%%%%%%%%%%%%%%%%%%%%%%%%%

\bibliographystyle{unsrt}  
\bibliography{references}  

%%%%%%%%%%%%%%%%%%%%%%%%%%%%%%%%%%%%%%%%%%%%%%%%%%%%%%%%%%%%%%%%%%%%%%%%%%%%%%%%%%%%%%%%

\newpage

\renewcommand{\appendixpagename}{\Large{Appendix}}
\appendixpage

\begin{appendices} \label{appendix}
\crefalias{section}{appendix}
%\numberwithin{equation}{section}

\input{sec6}

\end{appendices}

%%%%%%%%%%%%%%%%%%%%%%%%%%%%%%%%%%%%%%%%%%%%%%%%%%%%%%%%%%%%%%%%%

\end{document}

%% file: sec1.tex
% !TeX root = ms.tex

\section{Introduction} \label{sec:intro}

%parallelisation and why we use time-parallel methods.
In its most basic form, parallel computing is the process by which an algorithm is partitioned into a number of sub-problems that can be solved simultaneously without prior knowledge of each other. More widespread parallelism is becoming increasingly necessary in many different fields to reduce the computational burden and thus overcome the physical limitations (i.e. space, power usage, clock speeds, cooling, and financial costs) arising on machine hardware. Complex models in science often involve solving large systems of ordinary or partial differential equations (ODEs or PDEs) which, in the case of spatio-temporal PDEs, can be parallelised using existing domain decomposition methods. We refer to \cite{Toselli2005} for an overview. Although very efficient for high dimensional systems, these methods are reaching scale-up limits, and wallclock speeds often bottleneck due to the temporal integration. For instance, modern algorithms used to simulate Edge Localised Modes in turbulent fusion plasmas can take anywhere between 100-200 days to integrate over a time interval of just one second \cite{Samaddar2019}.

%These methods usually involve partitioning the spatial domain to solve smaller sub-problems in parallel before the solutions are pieced back together and the next time step is considered sequentially.

%Converting existing sequential algorithms into their equivalent parallel counterparts is challenging and efficiency gains will inevitably plateau. However, it is the dominant method for relieving the aforementioned strains on computer hardware.

%what is time parallelism and when did it begin, what methods are out there?
These sequential bottlenecks in time have motivated the research and development of \textit{time-parallel} integration methods for systems of ODEs and PDEs over the last 55 years or so (see \cite{Burrage1995,Gander2015,Ong2020} for reviews). These methods provide a way to integrate initial value problems (IVPs) over long time intervals where solutions would be unobtainable using existing sequential algorithms. One approach to integrate in a non-sequential manner, similar to spatial parallelisation, is to discretise the time interval into $N$ sub-intervals upon which $N$ sub-problems are solved in parallel using existing numerical methods. The solution at a given time step, however, depends upon the solution at the previous step(s). This creates a problem prior to the parallel integration as $N-1$ of the $N$ initial values (from which to begin integration in each sub-interval) are unknown \textit{a priori}. Existing algorithms that attempt to locate these $N-1$ initial values by direct or iterative means have been collectively referred to as multiple shooting methods \cite{Bellen1989,Chartier1993,Lions2001,Nievergelt1964a,Saha1997}.

%roughly what is parareal and what has been done with it?
Our focus is on one multiple shooting method in particular: \textit{parareal} \cite{Lions2001}. This algorithm uses a combination of coarse and fine grained (in time) numerical integrators to provide parallel speedup, and has been tested successfully on problems spanning molecular \cite{Baffico2002} and fluid dynamics \cite{Fischer2005,Garrido2005,Trindade2006}, to stochastic differential equations \cite{Engblom2009,Legoll2020} and fusion plasma simulations \cite{Samaddar2019,Samaddar2010}. The popularity of parareal is due to its relatively straightforward implementation and demonstrable effectiveness in speeding up the integration of IVPs. The algorithm iteratively locates a numerical solution at each sub-interval using a predictor-corrector scheme, derived by discretising the Newton-Raphson method. The algorithm stops once a tolerance is met with convergence guaranteed under certain mathematical conditions -- more detail on the algorithm is provided in \Cref{sec:parareal}. Much work has gone into analysing the numerical convergence of this method \cite{Bal2005,Bal2002,Gander2007a,Maday2002}, combining it with spatial decomposition techniques \cite{Maday2005}, and developing variants that utilise idle processors \cite{Elwasif2011} and adaptive time-stepping \cite{Maday2020}. However, given fixed fine and coarse integrators, the method itself is strictly deterministic and little work has gone into trying to reduce the number of iterations $k_d \in \{1,\dots,N \}$ that parareal takes to solve a particular IVP. Our primary focus is on reducing $k_d$ and henceforth we refer to it as the `convergence rate'. In scenarios where parareal struggles to converge in $k_d \ll N$ iterations, reducing $k_d$ by even a few iterations, say to $k_s < k_d$, can lead to significant parallel gains (roughly a factor $k_d/k_s$), enabling faster numerical integration.

%put some of this in above
%The parareal algorithm is a time-parallel method for solving systems of ODEs in a multiple shooting-type manner. Much like some of the domain decomposition methods in \cite{Toselli2005}, the idea behind this particular algorithm is to partition the time domain into smaller, non-overlapping sub-intervals upon which a set of sub-problems can be independently solved in parallel. For ODE problems however, the inherent sequential nature of time demands that initial values must be known at the start of each sub-interval before parallel integration can begin. Parareal provides a robust way of initially locating these initial values, converging to an accurate solution.

%An iterative time-parallel integrator, it was first introduced in \cite{Lions2001} and has grown increasingly popular since its inception due to its flexibility and ease of implementation to a number of applied mathematical problems.

%what do we do? what is stochastic parareal?
The purpose of this work is to extend the parareal algorithm using probabilistic methods to converge in fewer than $k_d$ iterations for a given system of ODEs. To achieve this, we introduce a \textit{stochastic} parareal algorithm. Instead of integrating forward in time from a single initial value (given by the predictor-corrector) in each sub-interval, a pre-specified probability distribution (see sampling rules in \Cref{subsec:sampling}) is used to generate $M$ candidate initial values that are \textit{simultaneously} integrated forward in time in parallel. At each sub-interval, one optimal initial value (from the $M$ samples) is selected sequentially such that the most continuous trajectory is chosen across consecutive sub-intervals. These initial values are then fed directly into the predictor-corrector -- the idea being that this stochastically generated set of values provides a better guess to the solution than those found purely deterministically. For example, suppose running the predictor-corrector with initial values $\bm{x}^{(0)}$ yields convergence in $k_d$ iterations, generating the sequence of solutions $\{\bm{x}^{(0)},\bm{x}^{(1)},\dots,\bm{x}^{(k_d)} \}$. Instead of starting with $\bm{x}^{(0)}$, suppose we sample initial values from a probability distribution and choose some ``better" starting point which is close to, say, $\bm{x}^{(i)}$ for some $i \in \{1,\dots,k_d - 1 \}$. Then the sequence generated by the predictor-corrector would instead be approximately $\{\bm{x}^{(i)}, \bm{x}^{(i+1)}, \dots, \bm{x}^{(k_d)} \}$, converging in $k_d - i$ iterations. Therefore, given a fixed number of samples $M$, the stochastic parareal algorithm will converge in fewer than $k_d$ iterations with some non-zero probability.

% Recall the predictor-corrector (\ref{eq:pred_correc}) is a discrete approximation of the Newton-Raphson method for finding the roots of the nonlinear system (\ref{eq:exact}). Given a sufficiently accurate initial guess $\mathbf{U}^0_n$, $\forall n \in \{0,\dots,N\}$, $\mathcal{P}$ locates the sequence of increasingly accurate solutions $\{\mathbf{U}^0_n,\mathbf{U}^1_n,\dots,\mathbf{U}^{k_d}_n \}$. If instead we provide (\ref{eq:pred_correc}) with a "better" initial guess, say some randomly selected $\hat{\bm{\alpha}}^k_n$, $\forall n \in \{0,\dots,N\}$, during each iteration $k \geq 1$, then (\ref{eq:pred_correc}) will generate the sequence $\{\mathbf{U}^0_n,\hat{\bm{\alpha}}^1_n,\dots,\hat{\bm{\alpha}}^{k_s}_n \}$ and converge in $k_s \leq k_d$ iterations with some probability.

%Probabilistic numerics is a growing field and combining parareal with stochastic methods in order to achieve faster convergence is not something that has been considered before.

This idea of propagating multiple initial values in each sub-interval in parallel is not new. In the first major work proposing the time-parallel integration of IVPs, Nievergelt \cite{Nievergelt1964a} discussed choosing $M$ initial values deterministically. The method for determining the solution from this ensemble of trajectories was to combine two of the $M$ samples in each sub-interval using an interpolation coefficient determined from the preceding sub-interval. Whilst an excellent first incursion into the field, for nonlinear problems this direct method suffered from interpolation errors and questions remained on how to efficiently scale this up to systems of ODEs. We generalise the original idea of Nievergelt by combining it with the parareal algorithm, generating the $M$ initial values at each sub-interval using probability distributions based on information known about the solutions at the current iteration.

%Outline of the paper
The paper is organised as follows. In \Cref{sec:parareal}, we recall the parareal algorithm and its properties from a multiple shooting viewpoint, including a known modification that will enable us to extend the algorithm to incorporate stochastic sampling. In \Cref{sec:stochasticparareal}, we introduce the stochastic parareal algorithm. Following an explanation of its key features, we elucidate how a variety of different sampling rules can be flexibly interchanged in order to carry out the sampling. In \Cref{sec:results}, we conduct numerical experiments to illustrate the performance of stochastic parareal against its deterministic counterpart using the different sampling rules. Findings are presented for three ODE systems of increasing complexity, with two additional examples given in the Appendix. Results indicate that for sufficiently many samples, stochastic parareal almost certainly beats the convergence rate of the deterministic algorithm and generates (stochastic) solutions of comparable numerical accuracy. For the multivariate ODEs, results show that performance is improved by generating correlated, as opposed to uncorrelated, random samples. In \Cref{sec:conclusions}, conclusions are drawn and avenues for future research are discussed.

%% file: sec2.tex
% !TeX root = ms.tex

\section{The parareal algorithm}  \label{sec:parareal}
\begin{figure}[b!]
         \centering
         \includegraphics[width=0.99\textwidth]{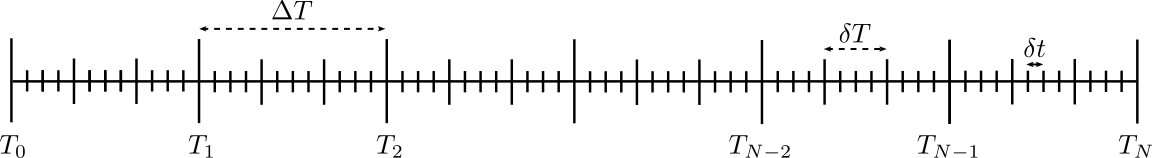}
         \caption{Schematic of the time domain decomposition. Three levels of temporal discretisation are shown: sub-intervals (size $\Delta T$), coarse intervals (size $\delta T$), and fine intervals (size $\delta t$). Note how the discretisations align with one another and that $\delta t \ll \delta T \leq \Delta T$ in our implementation.}
         \label{fig:time_decomp}
\end{figure}

Following previously outlined descriptions of the parareal algorithm  \cite{Gander2007a,Lions2001}, henceforth referred to as `parareal' or $\mathcal{P}$, consider a system of $d \in \mathbb{N}$ ODEs 
\begin{equation} \label{eq:ODE}
    \frac{d \mathbf{u}}{dt} = \mathbf{f} \big( \mathbf{u}(t),t \big)  \quad \text{on} \quad t \in [T_0, T_N], \quad \text{with} \quad \mathbf{u}(T_0) = \mathbf{u}^0,
\end{equation}    
where $\mathbf{f}: \mathbb{R}^d \times [T_0, T_N] \rightarrow \mathbb{R}^d $ is a smooth multivariate (possibly nonlinear) function, $\mathbf{u}: [T_0, T_N] \rightarrow \mathbb{R}^d$ the time-dependent vector solution, and $\mathbf{u}^0 \in \mathbb{R}^d$ the initial values at $T_0$. Decompose the time domain into $N$ sub-intervals such that $[T_0, T_N] = [T_0, T_1] \cup \dots \cup [T_{N-1}, T_N]$, with each sub-interval taking fixed length $\Delta T := T_n - T_{n-1}$ for $n = 1,\dots,N$ (see \Cref{fig:time_decomp} for a schematic).

Having partitioned the time domain, $N$ smaller sub-problems
\begin{equation} \label{eq:subprobs}
    \frac{d \mathbf{u}_n}{dt} = \mathbf{f} \big( \mathbf{u}_n (t | \mathbf{U}_n ),t \big) \quad \text{on} \quad t \in [T_n, T_{n+1}], \quad \text{with} \quad \mathbf{u}_n (T_n) = \mathbf{U}_n ,
\end{equation}
for $n = 0,\dots,N-1$ can theoretically be solved in parallel on $N$ processors, henceforth denoted $P_1,\dots,P_N$. Each solution $\mathbf{u}_n(t | \mathbf{U}_n)$ is defined over $[T_n,T_{n+1}]$ \textit{given} the initial values $\mathbf{U}_n \in \mathbb{R}^d$ at $t = T_n$. Note however that only the initial value $\mathbf{U}_0 = \mathbf{u}_0(T_0) = \mathbf{u}^0$ is known, whereas the rest ($\mathbf{U}_n$ for $n \geq 1$) need to be determined before \eqref{eq:subprobs} can be solved in parallel. These initial values must satisfy the continuity conditions
\begin{equation} \label{eq:exact}
    \mathbf{U}_0 = \mathbf{u}^0 \quad \text{and} \quad \mathbf{U}_n = \mathbf{u}_{n-1}(T_n | \mathbf{U}_{n-1}) \quad \text{for} \quad n = 1,\dots,N,
\end{equation}
which form a (nonlinear) system of $N+1$ equations that ensure solutions match at each $T_n \ \forall n \geq 1$. System \eqref{eq:exact} is solved for $\mathbf{U}_n$ using the Newton-Raphson method to form the iterative system
\begin{subequations} \label{eq:iterative}
\begin{align} 
    \mathbf{U}^{k+1}_0 &= \mathbf{u}^0, \label{eq:iterative_a} \\ 
    \mathbf{U}^{k+1}_n &= \mathbf{u}_{n-1} (T_n | \mathbf{U}^{k}_{n-1}) + \frac{\partial \mathbf{u}_{n-1}}{\partial \mathbf{U}_{n-1}} (T_n | \mathbf{U}^{k}_{n-1}) \big[\mathbf{U}^{k+1}_{n-1} - \mathbf{U}^{k}_{n-1} \big] , \label{eq:iterative_b}
\end{align}
\end{subequations}
for $n = 1,\dots,N$, where $k = 0,1,2,\dots$ is the iteration number. This system contains the unknown solutions $\mathbf{u}_n$ and their partial derivatives, which even if known, would be computationally expensive to calculate. 
\begin{figure}[t]
         \centering
         \includegraphics[width=0.99\textwidth]{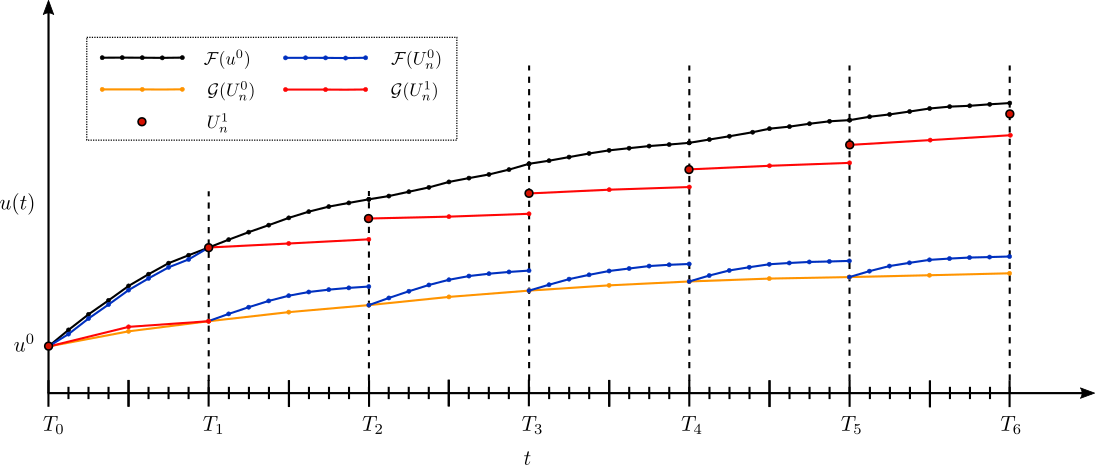}
         \caption{First iteration of the parareal algorithm to numerically evaluate the solution of a scalar ODE, either available or obtained via the fine solver (black line). The first runs of $\mathcal{G}$ and $\mathcal{F}$ are given in yellow and blue respectively; and the second run of $\mathcal{G}$ in red. The red dots represent the solution after applying the predictor-corrector \eqref{eq:pred_correc}.}
         \label{fig:parareal}
\end{figure}

To solve \eqref{eq:iterative} and evaluate $\mathbf{u}_n$ at discrete times, $\mathcal{P}$ utilises two numerical integrators. The first is a numerically fast \textit{coarse integrator} $\mathcal{G}$ that integrates over $[T_n,T_{n+1}]$ using initial values $\mathbf{U}_n$. The second is a \textit{fine integrator} $\mathcal{F}$ that runs much slower compared to $\mathcal{G}$ but has much greater numerical accuracy, chosen such that it integrates over the same interval $[T_n,T_{n+1}]$ with initial values $\mathbf{U}_n$. In our implementation, the difference between fast and slow integration is guaranteed by setting the time steps for $\mathcal{G}$ and $\mathcal{F}$ as $\delta T$ and $\delta t$ respectively with $\delta t \ll \delta T$ (recall \Cref{fig:time_decomp}). The key principle is that, if used to integrate \eqref{eq:ODE} over $[T_0,T_N]$ serially, $\mathcal{F}$ would  take an infeasible amount of computational time, hence the need to use $\mathcal{P}$ in the first place. Therefore $\mathcal{G}$ is permitted to run serially across multiple sub-intervals rapidly, whereas the slower solver $\mathcal{F}$ is only permitted to run in parallel on sub-intervals. This is a strict requirement for running $\mathcal{P}$, else numerical speedup will not be achieved.

Given that \eqref{eq:iterative_a} is known \textit{a priori} for all $k$, Lions et al. \cite{Lions2001} approximate the first term in equation \eqref{eq:iterative_b} using the fine solver $\mathcal{F}(\mathbf{U}^{k}_{n-1})$ and the second term by a coarse finite difference of the derivative using $\mathcal{G}(\mathbf{U}^{k+1}_{n-1}) - \mathcal{G}(\mathbf{U}^{k}_{n-1})$. One could instead approximate the derivative using a fine finite difference $\mathcal{F}(\mathbf{U}^{k+1}_{n-1}) - \mathcal{F}(\mathbf{U}^{k}_{n-1})$, however \eqref{eq:iterative_b} then becomes a sequential calculation using just $\mathcal{F}$ -- exactly what we wish to avoid. Assuming $\mathcal{G}$ meets the conditions required for \eqref{eq:iterative_b} to converge (see \Cref{subsec:convergence}), using the coarse approximation becomes reasonable and enables the parallel computation of \eqref{eq:iterative_b} \cite{Gander2007a}. The result is that an initial guess for the initial values $\mathbf{U}^{0}_n$ (found using $\mathcal{G}$) is improved at successive parareal iterations $k$ using the predictor-corrector update rule
\begin{equation} \label{eq:pred_correc}
        \mathbf{U}^{k+1}_n = \mathcal{G}(\mathbf{U}^{k+1}_{n-1}) + \mathcal{F}(\mathbf{U}^{k}_{n-1}) - \mathcal{G}(\mathbf{U}^{k}_{n-1}) \quad \text{for} \quad n = 1,\dots,N.
\end{equation}

%how it works
Pseudocode for an implementation of $\mathcal{P}$ is given in \Cref{alg:parareal} alongside a graphical description of the first iteration in \Cref{fig:parareal}. During the `zeroth' iteration, $\mathcal{G}$ is applied to $\mathbf{u}^0$, generating initial values $\mathbf{\hat{U}}^0_n \ \forall n \in \{1,\dots,N \}$ sequentially (lines 1-6). Immediately following, $\mathcal{F}$ is run in parallel from these initial values to generate a more accurate set of fine solution states $\mathbf{\tilde{U}}^0_n$ (lines 8-11). Next, $\mathcal{G}$ is run in the first time sub-interval, from the known initial value, to `predict' the solution $\mathbf{\hat{U}}^1_1$ at $T_1$. It is subsequently `corrected' using rule \eqref{eq:pred_correc}, and the predictor-corrector solution $\mathbf{U}^1_1$ is found at $T_1$. This prediction and correction process (lines 12-16) repeats sequentially until $\mathbf{U}^1_n$ is found $\forall n \in \{1,\dots,N \}$. After checking if the stopping criteria has been met (lines 17-21), the algorithm either starts the next iteration or stops.

%%% Parareal Algorithm
\begin{algorithm}
\caption{Parareal ($\mathcal{P}$)}
\label{alg:parareal}
\begin{algorithmic}[1]
\STATE{Set counters $k = I = 0$ and define $\mathbf{U}_n^k$, $\mathbf{\hat{U}}_n^k$ and $\mathbf{\tilde{U}}_n^k$ as the predictor-corrector, coarse, and fine solutions at the $n^{th}$ time step and $k^{th}$ iteration respectively (note that $\mathbf{U}_0^k = \mathbf{\hat{U}}_0^k = \mathbf{\tilde{U}}_0^k = \mathbf{u}^0 \ \forall k$).}
\STATE{// Calculate initial values at the start of each sub-interval $T_n$ by running $\mathcal{G}$ serially on processor $P_1$.}
\FOR{$n = 1$ \TO $N$} 
\STATE{$\mathbf{\hat{U}}_n^0 = \mathcal{G}(\mathbf{\hat{U}}_{n-1}^0)$} 
\STATE{$\mathbf{U}_n^0 = \mathbf{\hat{U}}_n^0$} 
\ENDFOR

\FOR{$k = 1$ \TO $N$}
\STATE{// Propagate the predictor-corrector solutions (from iteration $k-1$) on each sub-interval by running $\mathcal{F}$ in parallel on processors $P_{I+1},\dots,P_{N}$.}
\FOR{$n = I+1$ \TO $N$}
\STATE{$\mathbf{\tilde{U}}_n^{k-1} = \mathcal{F}(\mathbf{U}_{n-1}^{k-1})$} 
\ENDFOR
\STATE{// Propagate the predictor-corrector solution (at iteration $k$) with $\mathcal{G}$ on any available processor. Then correct this value using coarse and fine solutions obtained during iteration $k-1$ (this step \textit{cannot} be carried out in parallel).}
\FOR{$n = I+1$ \TO $N$}
\STATE{$\mathbf{\hat{U}}_n^k = \mathcal{G}(\mathbf{U}_{n-1}^{k})$} 
\STATE{$\mathbf{U}_{n}^{k} = \mathbf{\tilde{U}}_{n}^{k-1} + \mathbf{\hat{U}}_{n}^{k} - \mathbf{\hat{U}}_{n}^{k-1}$} 
\ENDFOR
\STATE{// Check whether the stopping criterion is met, saving all solutions up to time step $T_I$ before the next iteration. If tolerance is met for all time steps, the algorithm stops.}
\STATE{$I = \underset{n \in \{I+1,\dots,N\}}{\mathrm{max}}  \ \| \mathbf{U}^k_i - \mathbf{U}^{k-1}_i \|_{\infty} < \varepsilon \ \forall i < n$}
\IF{$I == N$}
\RETURN{$k$, $\mathbf{U}^k$}
\ENDIF
\ENDFOR
\end{algorithmic}
\end{algorithm}

\subsection{Stopping criteria and other properties} \label{subsec:convergence}
The algorithm is said to have converged up to time $T_I$ if
\begin{equation} \label{eq:convege}
    \| \mathbf{U}^k_n - \mathbf{U}^{k-1}_n \|_{\infty} < \varepsilon \quad \forall n \leq I,
\end{equation}
for some small fixed tolerance $\varepsilon$ -- noting that $\| \cdot \|_{\infty}$ denotes the usual infinity norm \cite{Maday2002,Samaddar2010}. Taking the relative, instead of absolute, errors in \eqref{eq:convege} would also be appropriate, however, in our numerical experiments the results (not reported) did not change. Once $I=N$, we say $\mathcal{P}$ has taken $k_d = k$ iterations to converge, yielding a solution with numerical accuracy of the order of the $\mathcal{F}$ solver. In its original formulation, $\mathcal{P}$ iteratively improves the solution across \textit{all} time steps, regardless of whether they have converged or not, up until the tolerance has been met for all $T_n$. The modified formulation presented here however only iterates on the solutions for the unconverged sub-intervals \cite{Elwasif2011,Samaddar2010}. This has no effect on the convergence rate $k_d$ and becomes especially important once we introduce the stochastic modifications in \Cref{sec:stochasticparareal}.

%(along with the solution given by $\mathcal{F}$ over $[T_0,T_N]$ which we are trying to locate)

%speedup notes
It should be clear that \textit{at least} one sub-interval will converge during each iteration $k$, as $\mathcal{F}$ will be run directly from a converged initial value. Therefore it will take at most $k_d=N$ iterations for $\mathcal{P}$ to converge, equivalent to running $\mathcal{F}$ over the $N$ sub-intervals serially (yielding zero parallel speedup). To achieve significant parallel speedup, multiple sub-intervals need to converge during an iteration so that $k_d \ll N$. Assuming $\mathcal{G}$ takes a negligible amount of time to run compared to the parallel components (along with any other serial computations), an approximate upper bound on the speedup achieved by $\mathcal{P}$ is $N/k_d$.

%discuss picking F and G
One challenge in attempting to minimise $k_d$, hence the overall runtime of $\mathcal{P}$, is to identify optimal solvers $\mathcal{F}$ and $\mathcal{G}$ for implementation. Whereas $\mathcal{F}$ is assumed to have high-accuracy and be computationally expensive to run, $\mathcal{G}$ must be chosen such that it runs significantly faster (usually orders of magnitude) than $\mathcal{F}$ whilst being sufficiently numerically accurate to converge in as few iterations as possible. Usually $\mathcal{G}$ is chosen such that its solutions have lower numerical accuracy, coarser temporal resolution, and/or reduced physics/coarser spatial resolution (when solving PDEs). There is currently no rigorous method for choosing the two solvers, however they should be chosen such that the ratio between the time taken to run $\mathcal{F}$ over $[T_{n},T_{n+1}]$ and the time taken to run $\mathcal{G}$ over the same interval is large.

In this paper, we fix both $\mathcal{F}$ and $\mathcal{G}$ as explicit\footnote{When solving the stiff ODE systems in \Cref{sec:results}, the initial values at $t=0$ are chosen such that explicit methods work. For greater numerical stability, one could alternatively use implicit solvers at extra computational cost.} fourth-order Runge-Kutta methods (henceforth RK4) for ease of implementation. More detailed analysis on the mathematical conditions required for $\mathcal{P}$ to converge, i.e. for the numerical solution $\mathbf{U}_n$ to approach the fine solution, can be found in \cite{Gander2008,Gander2007a,Lions2001,Maday2005}.

%% file: sec3.tex
% !TeX root = ms.tex

\section{A stochastic parareal algorithm} \label{sec:stochasticparareal}
In this section, we introduce the stochastic parareal algorithm, an extension of parareal that incorporates randomness and utilises its well-studied deterministic convergence properties to locate a solution in $k_s < k_d$ iterations. A summary of additional notation required to describe the algorithm, henceforth referred to as $\mathcal{P}_s$, is provided in \Cref{table:1}.

The idea behind $\mathcal{P}_s$ is to sample $M$ vectors of initial values $\bm{\alpha}^k_{n,1},\dots,\bm{\alpha}^k_{n,M}$, at each unconverged sub-interval $T_n$, in the neighbourhood of the current predictor-corrector solution $\mathbf{U}^k_n$ from a given probability distribution (with sufficiently broad support) and propagate them all in parallel using $\mathcal{F}$. Given a sufficient number of samples is taken, one will be closer (in the Euclidean sense) to the true root that equation \eqref{eq:pred_correc} is converging toward. Among them, we select an optimal $\hat{\bm{\alpha}}^k_n$ by identifying which samples generate the most continuous trajectory, at the fine resolution, in phase space across $[T_0,T_N]$. Therefore, at each iteration, we stochastically ``jump" toward more accurate initial values and then feed them into the predictor-corrector \eqref{eq:pred_correc}. For increasing values of $M$, the convergence rate $k_s$ will decrease, satisfying $k_s < k_d$ with probability one -- shown numerically in \Cref{sec:results}.

\begin{table}[tbhp]
{\footnotesize
\caption{Additional notation used to describe the stochastic parareal algorithm.} \label{table:1}
\begin{center}
\begin{tabular}{|c|p{0.8\linewidth}|} \hline \hline
\textbf{Notation} & \textbf{Description} \\ \hline
$n$    & Index of the time sub-interval $T_n$, $ n = 0,\dots,N $. \\ \hline
$k$    & Iteration number of $\mathcal{P}$ and $\mathcal{P}_s$, $k = 1,\dots,N $. \\ \hline
$k_s$  & Total iterations taken for stochastic parareal to stop and return a solution, $k_s \in \{1,\dots,N \} $. \\ \hline
$M$    & Number of random samples taken at each $T_n$. \\ \hline
$\Phi^{k-1}_{n}$ & The $d$-dimensional probability distribution used to sample initial values at time $T_n$ and iteration $k$. \\ \hline
$\bm{\alpha}^{k-1}_{n,m}$  & The $m^{\text{th}}$ $d$-dimensional sample from distribution $\Phi^{k-1}_{n}$ at time $T_n$ and iteration $k$, $m = 1,\dots,M$. \\ \hline
$\hat{\bm{\alpha}}^{k-1}_n$    & The selected $d$-dimensional sample from $\bm{\alpha}^{k-1}_{n,1},\dots,\bm{\alpha}^{k-1}_{n,M}$ at time $T_n$ and iteration $k$. \\ \hline
$\bm{\mu}^{k-1}_n$ & The $d$-dimensional vector of marginal means at time $T_n$ and iteration $k$. \\ \hline
$\bm{\sigma}^{k-1}_n$ & The $d$-dimensional vector of marginal standard deviations at time $T_n$ and iteration $k$. \\ \hline
$\rho^{k-1}_{n_{i,j}}$  & Pairwise correlations coefficients between $i^{\text{th}}$ and $j^{\text{th}}$ components of the $M$ fine resolution propagated samples at time  $T_n$ and iteration $k$. \\ \hline
$\mathbf{R}_n^{k-1}$ & The symmetric positive semi-definite $d \times d$ correlation matrix with elements $\rho^{k-1}_{n_{i,j}}$ at time $T_n$ and iteration $k$. \\ \hline
$\mathbb{I}$ & The $d \times d$ identity matrix. \\ \hline
$\bm{\Sigma}_n^{k-1}$ & The symmetric positive semi-definite $d \times d$ covariance matrix at time $T_n$ and iteration $k$. \\
\hline \hline
\end{tabular}
\end{center}
}
\end{table}

\subsection{The algorithm} \label{subsec:stochastic_parareal}
Suppose again we aim to solve system \eqref{eq:ODE}, adopting the same conditions, properties, and notation as discussed in \Cref{sec:parareal}. $\mathcal{P}_s$ follows the first iteration ($k=1$) of $\mathcal{P}$ identically -- see line 1 of \Cref{alg:stochparareal}. This is because information about initial values at the different temporal resolutions (i.e. results from $\mathcal{F}$ and $\mathcal{G}$) are required to construct the appropriate probability distributions for sampling. Following the convergence check, we assume (for the purposes of explaining the stochastic iterations) that only the first sub-interval $[T_0,T_1]$ converged during $k=1$, leaving $N-1$ unconverged sub-intervals. At this point we know the most up-to-date predictor-corrector solutions $\mathbf{U}_n^1$ $\forall n \in \{1,\dots,N \}$ and the stochastic iterations can begin (henceforth $k=2$).
\begin{figure}[t]
         \centering
         \includegraphics[width=0.99\textwidth]{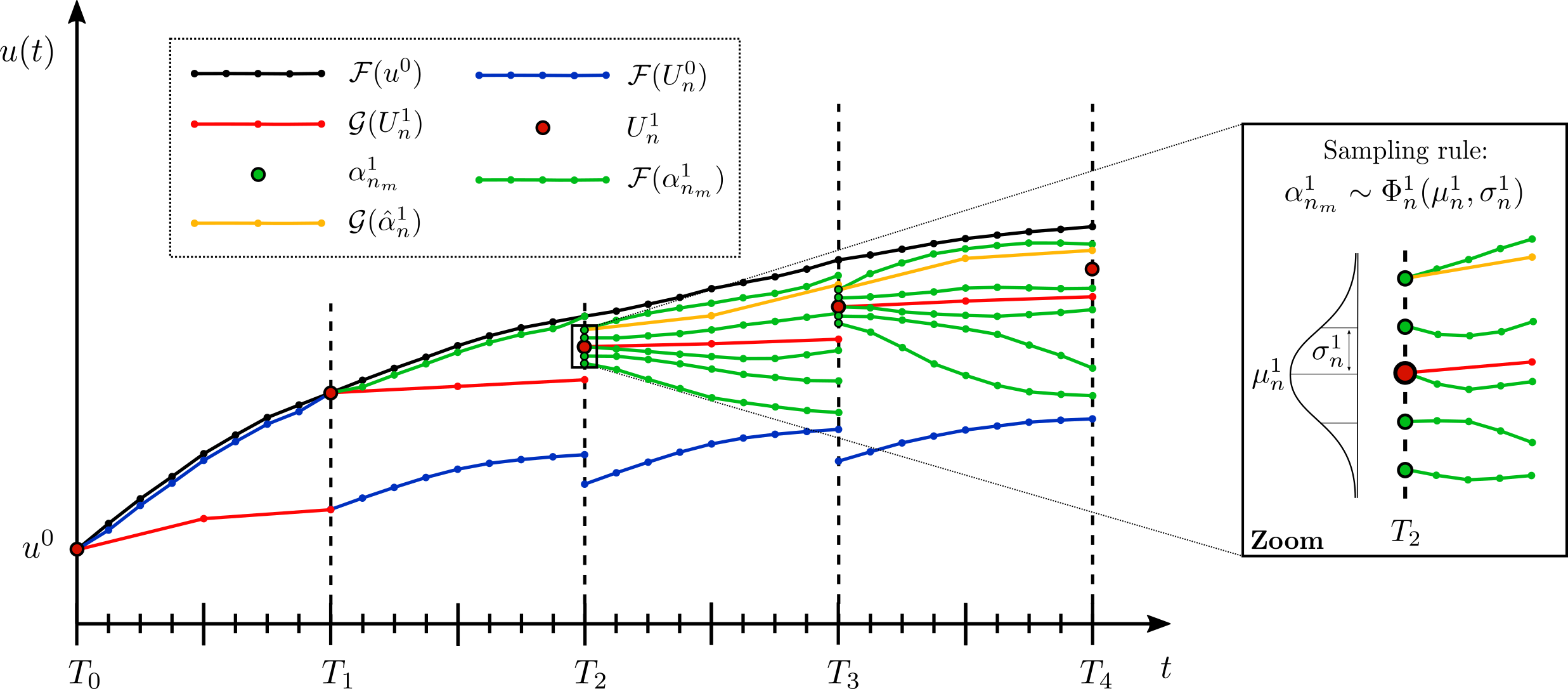}
         \caption{An illustration of the sampling and propagation process within $\mathcal{P}_s$ following iteration $k=1$. The fine solution is given in black, the $k=0$ fine solutions in blue, the $k=1$ coarse solutions in red, and the $k=1$ predictor-corrector solutions as red dots. With $M=5$, four samples $\alpha^{1}_{n,m}$ (green dots) are taken at $T_2$ and $T_3$ from distributions with means $U^1_2$ and $U^1_3$, and some finite standard deviations respectively. These values, along with $U^1_2$ and $U^1_3$ themselves, are propagated in parallel forward in time using $\mathcal{F}$ (green lines). The optimally chosen samples $\hat{\alpha}^1_n$ (refer to text for how these are chosen) are then propagated forward in time using $\mathcal{G}$ (yellow lines).}
         \label{fig:stoch_parareal}
\end{figure}

%the correlations and sampling
At any unconverged $T_n$ ($n>1$), we sample $M$ vectors of initial values, denoted $\bm{\alpha}^{k-1}_{n,m}$ for $m = 1,\dots,M$. The first sample is fixed as the predictor-corrector $\mathbf{U}^{k-1}_n$, to ensure that $\mathcal{P}_s = \mathcal{P}$ when $M=1$. The other $M-1$ initial values are sampled from a pre-specified $d$-dimensional probability distribution $\Phi^{k-1}_n$ with finite marginal means $\bm{\mu}_n^{k-1} = (\mu_{n_1}^{k-1},\dots,\mu_{n_d}^{k-1})^{\text{T}}$, marginal standard deviations $\bm{\sigma}_{n}^{k-1} = (\sigma_{n_1}^{k-1},\dots,\sigma_{n_d}^{k-1})^{\text{T}}$, and correlation structure given by the matrix $\mathbf{R}^{k-1}_n$. These quantities depend upon the information available at iteration $k-1$, i.e. a combination of $\mathbf{U}^{k-1}_n$, $\mathcal{F}(\mathbf{U}^{k-1}_n)$, $\mathcal{G}(\mathbf{U}^{k-1}_n)$ and $\mathcal{G}(\mathbf{U}^{k-2}_n)$, see \Cref{subsec:sampling}. The correlation matrix $\mathbf{R}^{k-1}_n$ is introduced to take into account the dependence between components of the ODE system (lines 3-9). The elements of $\mathbf{R}^{k-1}_n$, for $k \geq 3$, are defined using the Pearson correlation coefficient
\begin{equation} \label{eq:corrs}
    \rho_{n_{i,j}}^{k-1} = \frac{\sum_{m=1}^M (x_m^{(i)} - \bar{x}^{(i)}) (x_m^{(j)} - \bar{x}^{(j)})}{\sqrt{\sum_{m=1}^M (x_m^{(i)} - \bar{x}^{(i)})^2} \sqrt{\sum_{m=1}^M (x_m^{(j)} - \bar{x}^{(j)})^2 }} , \quad i,j \in \{1,\dots,d\},
\end{equation}
where 
\begin{align}
    x_m^{(i)} = \mathcal{F}(\bm{\alpha}^{k-2}_{{n-1},m})(i), \qquad \bar{x}^{(i)} = \frac{1}{M} \sum_{m=1}^M x_m^{(i)},
\end{align}
and $\mathcal{F}(\bm{\alpha}^{k-2}_{{n-1},m})(i)$ denotes the $i$th element of $\mathcal{F}(\bm{\alpha}^{k-2}_{{n-1},m})$.
The coefficients $ \rho_{n_{i,j}}^{k-1}$ in \eqref{eq:corrs} are essentially the estimated pairwise correlation coefficients of the $M$ $d$-dimensional fine resolution propagations of the sampled initial values at $T_n$ from the previous iteration, i.e.\ $\mathcal{F}(\bm{\alpha}^{k-2}_{{n-1},1}),\ldots,\mathcal{F}(\bm{\alpha}^{k-2}_{{n-1},M})$. 
Note that other types of linear correlation coefficient could also be chosen. 
Since each $\mathcal{F}(\bm{\alpha}^{k-2}_{{n-1},m})$ is not available at iteration $k=2$, we set $\mathbf{R}_n^{k-1} = \mathbb{I}$ for $k=2$, i.e. we sample from a multivariate distribution with uncorrelated components.

Following this, the sampling and subsequent propagation using $\mathcal{F}$ can begin in parallel (lines 10-22). Given the solution between $[T_0,T_1]$ has converged, $\mathcal{F}$ will run from the converged initial value at $T_1$, with sampling starting from $T_2$ onward (see \Cref{fig:stoch_parareal}). \textit{All} sampled initial values are then propagated forward in parallel using $\mathcal{F}$, requiring at least $M(N-2)+1$ processors ($M$ samples times $N-2$ unconverged sub-intervals plus running $\mathcal{F}$ once in $[T_1,T_2]$).

Of the $M$ sampled initial values at each $T_n$ ($n>1$), only one is retained, denoted by $\bm{\hat{\alpha}}^{k-1}_n$, chosen such that it minimises the Euclidean distance between the fine solution and the sampled values (lines 23-28). To do this, start from the converged initial values at $T_2$ given by the fine solver: $\mathcal{F}(\mathbf{U}^{k-1}_1)$. Calculate the Euclidean distance between $\mathcal{F}(\mathbf{U}^{k-1}_1)$ and each of the $M$ samples $\bm{\alpha}^{k-1}_{2,1},\dots,\bm{\alpha}^{k-1}_{2,M}$. The sample minimising this distance is chosen as $\bm{\hat{\alpha}}^{k-1}_2$. Repeat for later $T_n$, minimising the distance between $\mathcal{F}(\hat{\bm{\alpha}}^{k-1}_{n-1})$ and one of the samples $\bm{\alpha}^{k-1}_{n,1},\dots,\bm{\alpha}^{k-1}_{n,M}$. This process must be run sequentially and relies on the modification to $\mathcal{P}$ made in \Cref{sec:parareal} -- that solutions are not altered once converged. Referring again to \Cref{fig:stoch_parareal}, the corresponding coarse trajectories of these optimally chosen samples $\hat{\bm{\alpha}}^{k-1}_n$ must also be calculated to carry out the predictor-corrector step (lines 29-32).

At this point, the set of initial values $\{ \bm{\hat{\alpha}}^{k-1}_2, \dots, \bm{\hat{\alpha}}^{k-1}_{N-1} \}$ has been selected from the ensemble of random samples, effectively replacing the previously found $\mathbf{U}^{k-1}_n$. The coarse and fine propagations of these values are now used in the predictor-corrector (lines 33-37) such that
\begin{equation} \label{eq:pred_correc_stoch}
        \mathbf{U}^{k}_n = \left\{
	\begin{array}{ll}
		\mathcal{G}(\mathbf{U}^{k}_{n-1}) + \mathcal{F}(\mathbf{U}^{k-1}_{n-1}) - \mathcal{G}(\mathbf{U}^{k-1}_{n-1}) & \mbox{for} \quad n = 2, \\
		\mathcal{G}(\mathbf{U}^{k}_{n-1}) + \mathcal{F}(\bm{\hat{\alpha}}^{k-1}_{n-1}) - \mathcal{G}(\bm{\hat{\alpha}}^{k-1}_{n-1}) & \mbox{for} \quad n = 3,\dots,N. 
	\end{array}
\right.
\end{equation}
Using the same stopping criteria \eqref{eq:convege} from $\mathcal{P}$ (lines 38-42), the algorithm either stops or runs another stochastic parareal iteration.

%%% Stochastic Parareal Algorithm
\begin{algorithm}
\caption{Stochastic parareal ($\mathcal{P}_s$)}
\label{alg:stochparareal}
\begin{algorithmic}[1]
\STATE\label{line3}{Run $\mathcal{P}$ until the end of iteration $k=1$.}
\FOR{$k = 2$ \TO $N$}
\STATE{// Calculate correlation matrices (only if $d > 1$).}
\STATE{$\mathbf{R}_n^{k-1} = \mathbb{I}$ $\forall n$}
\IF{$k \geq 3$}
\FOR{$n = I+1$ \TO $N-1$}
\STATE{Calculate $\mathbf{R}_n^{k-1}$ using $\mathcal{F}(\bm{\alpha}_{{n-1},1}^{k-2}),\dots,\mathcal{F}(\bm{\alpha}_{{n-1},M}^{k-2})$ (recall \eqref{eq:corrs}).}
\ENDFOR
\ENDIF
\STATE{// Initial value sampling and propagation. Both line 11 and the nested for loop below (lines 12--22) run in parallel on $P_1,\dots,P_{M(N-1-I)+1}$, i.e. all runs of $\mathcal{F}$ must be in parallel.}
\STATE{$\mathbf{\tilde{U}}^{k-1}_{I+1} = \mathcal{F}(\mathbf{U}^{k-1}_{I}) $ // propagate converged initial value at $T_I$ on $P_1$}
\FOR{$n = I+1$ \TO $N-1$}
\FOR{$m = 1$ \TO $M$}
\IF{$m == 1$}
\STATE{$\bm{\alpha}^{k-1}_{n,1} = \mathbf{U}^{k-1}_{n}$ // first `sample' is fixed to the predictor-corrector}
\STATE{$\mathbf{\tilde{U}}_{n+1,1} = \mathcal{F}(\bm{\alpha}_{n,1}^{k-1})$ // temporarily store propagated values}
\ELSE
\STATE{$\bm{\alpha}^{k-1}_{n,m} \sim \Phi^{k-1}_n$ // sample initial value randomly, see \Cref{subsec:sampling}}
\STATE{$\mathbf{\tilde{U}}_{n+1,m} = \mathcal{F}(\bm{\alpha}_{n,m}^{k-1})$}
\ENDIF
\ENDFOR
\ENDFOR
\STATE{// Select the most continuous fine trajectory from the ensemble sequentially.}
\FOR{$n = I+1$ \TO $N-1$} 
\STATE{$ J = \underset{j \in \{1,\dots,M\}}{\mathrm{arg min}} \ \| \bm{\alpha}^{k-1}_{n,j} - \mathbf{\tilde{U}}^{k-1}_{n} \|_2 $} 
\STATE{$\hat{\bm{\alpha}}_{n}^{k-1} = \bm{\alpha}^{k-1}_{n,J} $ // store optimal initial value} 
\STATE{$\mathbf{\tilde{U}}_{n+1}^{k-1} = \mathbf{\tilde{U}}_{n+1,J}$ // store most optimal fine trajectories}
\ENDFOR
\STATE{// Run $\mathcal{G}$ from the optimal samples (can run in parallel).}
\FOR{$n = I+1$ \TO $N-1$} 
\STATE{$\mathbf{\hat{U}}_{n+1}^{k-1} = \mathcal{G}(\hat{\bm{\alpha}}_{n}^{k-1})$} 
\ENDFOR
\STATE{// Predict and correct the initial values.}
\FOR{$n = I+1$ \TO $N$}
\STATE{$\mathbf{\hat{U}}_n^{k} = \mathcal{G}(\mathbf{U}_{n-1}^{k})$} 
\STATE{$\mathbf{U}_n^{k} = \mathbf{\tilde{U}}_{n}^{k-1} + \mathbf{\hat{U}}_{n}^{k} - \mathbf{\hat{U}}_{n}^{k-1}$} 
\ENDFOR
\STATE{// Check whether the stopping criterion is met.}
\STATE{$I = \underset{n \in \{I+1,\dots,N\}}{\mathrm{max}} \ \| \mathbf{U}^k_i - \mathbf{U}^{k-1}_i \|_{\infty} < \varepsilon \ \forall i < n$}
\IF{$I == N$}
\RETURN{$k$, $\mathbf{U}^k$}
\ENDIF
\ENDFOR
\end{algorithmic}
\end{algorithm}

\subsection{Sampling rules} \label{subsec:sampling}
The probability distributions $\Phi^{k-1}_n$ incorporate different combinations of available information about the initial values at different temporal resolutions, i.e. the coarse, fine, and predictor-corrector initial values $\mathcal{G}(\mathbf{U}^{k-1}_{n-1})$, $\mathcal{G}(\mathbf{U}^{k-2}_{n-1})$, $\mathcal{F}(\mathbf{U}^{k-2}_{n-1})$, and $\mathbf{U}^{k-1}_{n}$ respectively. This information is used to define the marginal means and standard deviations in the `sampling rules' outlined in the following subsections. Note how the distributions vary with time $T_n$ and $k$, so that the accuracy of the distributions increase as the initial values in $\mathcal{P}_s$ are iteratively improved. Using Gaussian and copula distributions, we analyse the performance of different sampling rules within $\mathcal{P}_s$ in \Cref{sec:results}. This will give us a more comprehensive understanding on whether the choice of distribution family $\Phi^{k-1}_n$ or the parameters $\bm{\mu}_n^{k-1}$, $\bm{\sigma}_n^{k-1}$, and $\mathbf{R}^{k-1}_n$ have the greatest impact on the convergence rate $k_s$. In all tested cases, we observe that taking correlated samples significantly improves the performance of $\mathcal{P}_s$ compared to the independent setting.

\subsubsection{Multivariate Gaussian} \label{subsubsec:gaussians}
First, we consider perturbing the initial values using stochastic ``noise", i.e. considering errors compared to the true initial values to be normally distributed, a standard method for modelling uncertainty. The initial values $\bm{\alpha}^{k-1}_{n,m}$ are sampled from a multivariate Gaussian distribution $\mathcal{N}(\bm{\mu}_n^{k-1},\bm{\Sigma}_n^{k-1})$, where $(\bm{\Sigma}_n^{k-1})_{i,j} = \rho^{k-1}_{n_{i,j}} \sigma^{k-1}_{n_i} \sigma^{k-1}_{n_j}$ is the $d \times d$ positive semi-definite covariance matrix. As marginal means $\bm{\mu}_n^{k-1}$ we choose either the fine resolution values $\mathcal{F}(\mathbf{U}^{k-2}_{n-1})$ (prior to correction) or the predictor-corrector values $\mathbf{U}^{k-1}_{n}$. For the marginal standard deviations, we choose $\bm{\sigma}_n^{k-1} = |\mathcal{G}(\mathbf{U}^{k-1}_{n-1}) - \mathcal{G}(\mathbf{U}^{k-2}_{n-1})|$ as they are of the order of the corrections made by the predictor-corrector and each marginal decreases toward zero as the algorithm converges (as expected). For the correlation coefficients, $\rho^{k-1}_{n_{i,j}}$, we calculate the linear correlation between the $\mathcal{F}$ propagated samples using Pearson's method -- recall \eqref{eq:corrs}. Note that $| \cdot |$ denotes the component-wise absolute value. Testing revealed that alternative marginal standard deviations $| \mathbf{U}^{k-1}_{n} - \mathbf{U}^{k-2}_{n}|$ and $|\mathcal{F}(\mathbf{U}^{k-2}_{n-1}) - \mathcal{G}(\mathbf{U}^{k-2}_{n-1})|$ did not span sufficiently large distances around $\bm{\mu}^{k-1}_n$ in order for sampling to be efficient, i.e they required much higher sampling to perform as well as $\bm{\sigma}_n^{k-1} = |\mathcal{G}(\mathbf{U}^{k-1}_{n-1}) - \mathcal{G}(\mathbf{U}^{k-2}_{n-1})|$ (results not shown). The samples $\bm{\alpha}^{k-1}_{n,m} \sim \mathcal{N}(\bm{\mu}_n^{k-1},\bm{\Sigma}_n^{k-1})$ are taken according to the following sampling rules:
\begin{align*}
    \textbf{Rule 1}: \ \bm{\mu}_n^{k-1} &= \mathcal{F}(\mathbf{U}^{k-2}_{n-1}) \  \text{and} \ \bm{\sigma}_n^{k-1} = | \mathcal{G}(\mathbf{U}^{k-1}_{n-1}) - \mathcal{G}(\mathbf{U}^{k-2}_{n-1}) |. \\
    \textbf{Rule 2}: \ \bm{\mu}_n^{k-1} &= \mathbf{U}^{k-1}_{n} \  \text{and} \ \bm{\sigma}_n^{k-1} = | \mathcal{G}(\mathbf{U}^{k-1}_{n-1}) - \mathcal{G}(\mathbf{U}^{k-2}_{n-1}) |.
\end{align*}
Note that a linear combination of the rules or taking half the samples from each appears to work well, with performance similar to the individual rules themselves (results not shown).

\subsubsection{Multivariate copula}
%what is a copula?
Alternatively, we consider errors that may have a different (possibly non-Gaussian) dependency structure. We consider another multivariate distribution, known as a copula, with uniformly distributed marginals scaled such that they have the same value as the marginal means and standard deviations in \Cref{subsubsec:gaussians}.

A copula $\mathcal{C}: [0,1]^d \rightarrow [0,1]$ is a joint cumulative distribution function, of a $d$-dimensional random vector, with uniform marginal distributions \cite{Nelsen2006}. Sklar's theorem \cite{Sklar1959} states that any multivariate cumulative distribution function with continuous marginal distributions can be written in terms of $d$ uniform marginal distributions and a copula that describes the correlation structure between them. Whilst there are numerous families of copulas, we consider the symmetric $t$-copula, $\mathcal{C}^{t}$, the copula underlying the multivariate $t$-distribution, which depends on the parameter $\nu$ and the correlation matrix $\mathbf{R}_n^{k-1}$. We fix $\nu = 1$ so that samples have a higher probability of being drawn toward the edges of the $[0,1]^d$ hypercube, see \cite{Nelsen2006}. Note that $\nu \rightarrow \infty$ can be thought of as sampling from the Gaussian copula. 

%scaling the copula
Correlated samples $\bm{\chi} = (\chi_1,\dots,\chi_d)^{\text{T}} \sim \mathcal{C}^{t}$ that are generated in $[0,1]^d$ then need to be re-scaled such that each marginal is uniformly distributed in an interval $[x_i,y_i] \subset \mathbb{R}$, with mean $\mu_i$ and standard deviation $\sigma_i$ for $i \in \{1,\dots,d \}$. By definition, a marginal uniform distribution on $[x_i,y_i]$ has mean $(x_i + y_i)/2$ and variance $(y_i - x_i)^2 /12$ which we set to $\mu_i$ and $\sigma^2_i$ respectively. Solving these equations, we find that the desired marginals are uniform distributions on $[\mu_i - \sqrt{3} \sigma_i ,\mu_i + \sqrt{3} \sigma_i]$. Thus, setting $2 \sqrt{3} \sigma_i \chi_i + \mu_i - \sqrt{3} \sigma_i$ guarantees that the generated samples $\bm{\chi} \sim \mathcal{C}^t$ have the same marginal means $\mu_i$ and standard deviations $\sigma_i$ as the Gaussian distributions. This allows us to compare performance results in \Cref{sec:results}. The $t$-copula sampling rules (\textbf{Rule 3} and \textbf{Rule 4}) are therefore defined component-wise as $\bm{\alpha}_{n,m}^{k-1}(i) = 2 \sqrt{3} \sigma_i \chi_i + \mu_i - \sqrt{3} \sigma_i$, for $i \in \{1,\dots,d \}$, with $\bm{\chi} \sim \mathcal{C}^{t}$ and parameters $\mu_i$ and $\sigma_i$ chosen to be identical to Rule 1 and Rule 2 respectively.

Note that this copula construction can be extended to other families of copulas (e.g.\ the Gaussian copula) and to copulas with different marginal distributions (e.g.\ Gaussian, $t$, or logistic marginals). 
Therefore, sampling from the Gaussian multivariate distributions in \Cref{subsubsec:gaussians} can be considered a special case of these more general copulas.

%some remarks on convergence
\subsection{Convergence} \label{subsec:SP_convergence}
Independent simulations of $\mathcal{P}_s$ will return a numerical solution $\bm{U}^{k}_n$ which, along with $k_s$, vary stochastically -- since the optimal initial values chosen will vary between simulations. Given the randomness in these quantities, we can discuss the convergence of $\mathcal{P}_s$ in two ways. 

% \subsubsection{Minimising $k_s$}
Firstly, consider convergence in terms of minimising the random variable $k_s$ by studying $\mathbb{P}(k_s<k_d)$. Given there are no analytical results for $\mathcal{P}$ guaranteeing that $k_d<N$ for any given problem, proving that $\mathbb{P}(k_s<k_d) = 1$, or at least $\mathbb{E}(k_s) = \sum_{k=1}^N k \mathbb{P}(k_s = k) < k_d$ seems to be analytically intractable. 
However, we can qualitatively discuss $\mathbb{P}(k_s<k_d)$ and $\mathbb{E}(k_s)$ with respect to the number of samples $M$. Consider the following cases: 
\begin{itemize}
    \item $M = 1$ \\
    Running $\mathcal{P}_s$ is equivalent to running $\mathcal{P}$, hence the convergence of $\mathcal{P}_s$ follows from that of $\mathcal{P}$ and therefore $\mathbb{E}(k_s) = k_s = k_d$.
    
    \item $1 < M < \infty$ \\
    For finitely many samples, we estimate the discrete probability distributions $\mathbb{P}(k_s=k)$ and observe, in all numerical experiments (see \Cref{sec:results}), that $\mathbb{P}(k_s<k_d) \rightarrow 1$ for $M \approx 10$.
    Moreover, we observe that $\mathbb{P}(k_s<k_d)$ increases and $\mathbb{E}(k_s)$ decreases for increasing $M$ -- with $\mathbb{E}(k_s) < k_d$ for all values of $M$ tested.
    
    \item $M \rightarrow \infty$ \\
    If we were able to take infinitely many samples, $\mathcal{P}_s$ effectively samples every possible value in the support of $\Phi$, i.e. every initial value that has a non-zero probability of being sampled from $\Phi$. Therefore, if $\Phi$ has infinite support, e.g. the Gaussian distribution, all possible initial values in $\mathbb{R}^d$ are sampled and propagated, hence the fine solution will be recovered almost surely in $\mathbb{E}(k_s) = k_s = 2$ iterations. Note this is the smallest value $k_s$ can take to converge assuming convergence does not occur following the first iteration -- although the algorithm could be modified such that convergence occurs almost surely in $k_s=1$ iterations as $M \rightarrow \infty$. In \Cref{subsec:1Dnonlinear}, we illustrate this property numerically by taking a large number of samples for a single realisation of $\mathcal{P}_s$.
\end{itemize}
In the scenario that $\mathcal{P}_s$ converges in $k_s=N$ iterations (irrespective of the value of $M$), it will return the fine solution just as $\mathcal{P}$ does when $k_d=N$ (having propagated the exact initial value at $T_0$ sequentially $N$ times using $\mathcal{F}$).

Secondly, consider convergence in terms of the stochastic solution $U^k_n$ \eqref{eq:pred_correc_stoch} approaching (in the mean-square sense) the fine solution $U_n$ as $k$ increases. 
In the numerical experiments in \Cref{sec:results}, we did not observe a case where $\mathcal{P}_s$ fails to converge to the fine solution.
In fact, we observe tight confidence intervals on the numerical errors between $U^k_n$ and $U_n$ upon multiple realisations of $\mathcal{P}_s$ (see \Cref{fig:chartier_errors} and \Cref{fig:bruss_errors}).
We may expect $\mathcal{P}_s$ to fail to converge (i.e.\ solutions blow up) in cases in which $\mathcal{P}$ also fails -- typically this means that a more accurate coarse solver is required for both algorithms.
It should be noted, however, that because $\mathcal{P}_s$ samples $M$ initial values, only $M'$ out of the $M$ propagated solutions may blow up in each time sub-interval, due to the poor coarse solver, and so $\mathcal{P}_s$ could in fact locate a solution in cases where $\mathcal{P}$ cannot -- although this is not tested here.

%some remarks on computational complexity
\subsection{Computational complexity}
At each iteration, $\mathcal{P}_s$ runs the fine solver more frequently than $\mathcal{P}$, albeit still in parallel, and therefore requires a larger number of processors. The first iteration of $\mathcal{P}_s$ requires $N$ processors, however once sampling begins in $k \geq 2$, it requires at most $M(N-I-1)+1$ processors -- assuming $I$ sub-intervals converge during $k=1$. Whilst we assume processors are in abundance, this number scales directly with $M$ and so it is important to keep $M$ to a minimum if limited processing power is available in practice.

\begin{figure}[t!]
         \centering
         \includegraphics[width=0.99\textwidth]{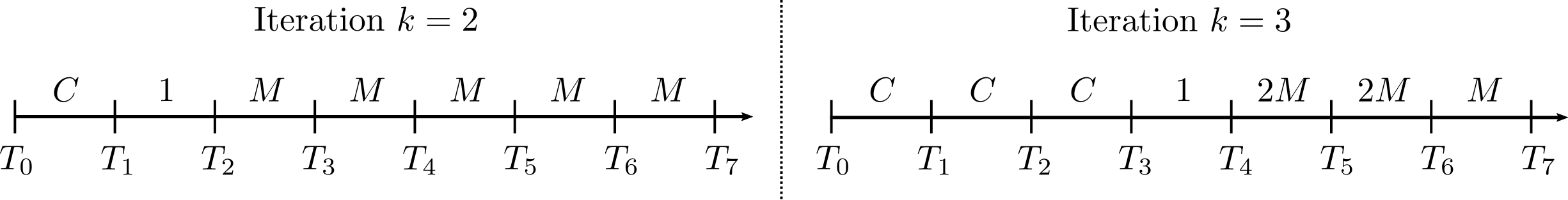}
         \caption{Illustration of a possible processor configuration if two sub-intervals were to converge between iterations two and three of $\mathcal{P}_s$. The letter $C$ denotes a converged sub-interval, the number $1$ denotes a sub-interval where we propagate the converged initial value from the preceding sub-interval using $\mathcal{F}$, and the letter $M$ denotes the number of samples taken in an unconverged sub-interval.}
         \label{fig:proc_config}
\end{figure}

As the stochastic iterations progress, the number of processors required, i.e.\ $M(N-I-1)+1$, decreases as the number of converged sub-intervals, $I$, increases -- leaving a growing number of processors idle. 
Each additional sub-interval that converges leaves $M$ idle processors that we can re-assign to do additional sampling and propagation.
We assign each set of $M$ idle processors to the earliest unconverged sub-interval with the least number of samples, ensuring all processors are working at all times to explore the solution space for the true initial values (see \Cref{fig:proc_config} for an illustration). 
We do not explicitly write the pseudocode for re-assigning the idle processors in \Cref{alg:stochparareal} (lines 10-22) to avoid additional complexity -- the process is, however, implemented in the numerical experiments in \Cref{sec:results}\footnote{To increase efficiency further we also attempted to store previously sampled and propagated ﬁne trajectories to use in future iterations of the algorithm, however, they did not improve performance. This was because only the most recent samples were ever chosen in each iteration (results not shown).}. 

In terms of timings, an iteration of $\mathcal{P}_s$ takes approximately the same wallclock time as one of $\mathcal{P}$. 
In this regard, we assume that the extra serial costs in $\mathcal{P}_s$, e.g. correlation estimation, selecting optimal samples, and the extra $\mathcal{G}$ runs, take negligible wallclock time when compared to a single run of $\mathcal{F}$. 
Therefore, the wallclock time for $\mathcal{P}_s$ will be lower than for $\mathcal{P}$ if $k_s < k_d$. 
This comes at a cost of requiring $\mathcal{O}(MN)$ processors rather than $N$ to solve the problem. 
However, it should be highlighted that if $\mathcal{P}_s$ converges in even one less iteration than $\mathcal{P}$, we avoid an extra run of $\mathcal{F}$ which may save a large amount of wallclock time.

%% file: sec4.tex
% !TeX root = ms.tex

\section{Numerical results} \label{sec:results}
%outline this section
In this section, we compare the numerical performance of $\mathcal{P}$ and $\mathcal{P}_s$ on systems of one, two, and three ODEs of varying complexity\footnote{All algorithms were coded in \texttt{MATLAB} and simulations were run using HPC facilities at the University of Warwick. Samples code for both $\mathcal{P}$ and $\mathcal{P}_s$ can be found in the public repository at \url{https://github.com/kpentland/StochasticParareal}.}. Both algorithms use RK4 methods to carry out integration, with $\mathcal{G}$ using time step $\delta T$ and $\mathcal{F}$ using time step $\delta t$, the latter at least 75 times smaller than the former. We quantify the performance of $\mathcal{P}_s$ by estimating the distributions of $k_s$ for each sampling rule and by measuring the accuracy of the stochastic solutions against those obtained serially with $\mathcal{F}$. Since a limited number of processors were available for these experiments, the results are based not on calculating wallclock runtimes but on comparing the convergence rates $k_d$ and $k_s$ -- which are independent of the number of processors used. Additional results for these test cases, as well as two further test problems, are given in the Appendix.

%solutions we could analyse:
% ODE solutions (phase planes or individual trajectories).
% could plot parareal successive errors for the 1D cases -- to show convergence.
% m vs k (for each rule on same plot).
% bar charts of the convergence distributions + a line giving the overall prob of beating parareal (for each rule -- then we can pick the best).
% errors vs. the fine solution for 250 runs (for a particular rule and given $M$) -- shows higher k are more accurate.
% varying ICs (u0 vs. k) -- can do, more difficult to then find m for each such that P(ks<kd)>0.95 -- doesn't need to be perfect, can choose some arbitrary Ndg that makes results look good.
%can plot scatter plots of the sampling (and which one gets selected) at each time sub-interval and iteration.
%alternatively: compare same rules for differing coarse steps (make note saying fine steps have no impact).

%%%%%%%%%%%%%%%%%%%%%%%%%%%%%%%%%%%%%%%%%%%%%%%%%%%%%%%%%%%%%%%%%
\subsection{Scalar nonlinear equation} \label{subsec:1Dnonlinear}
\begin{figure}[t!]
    \centering
    \subfloat[]{\includegraphics[width=0.49\textwidth]{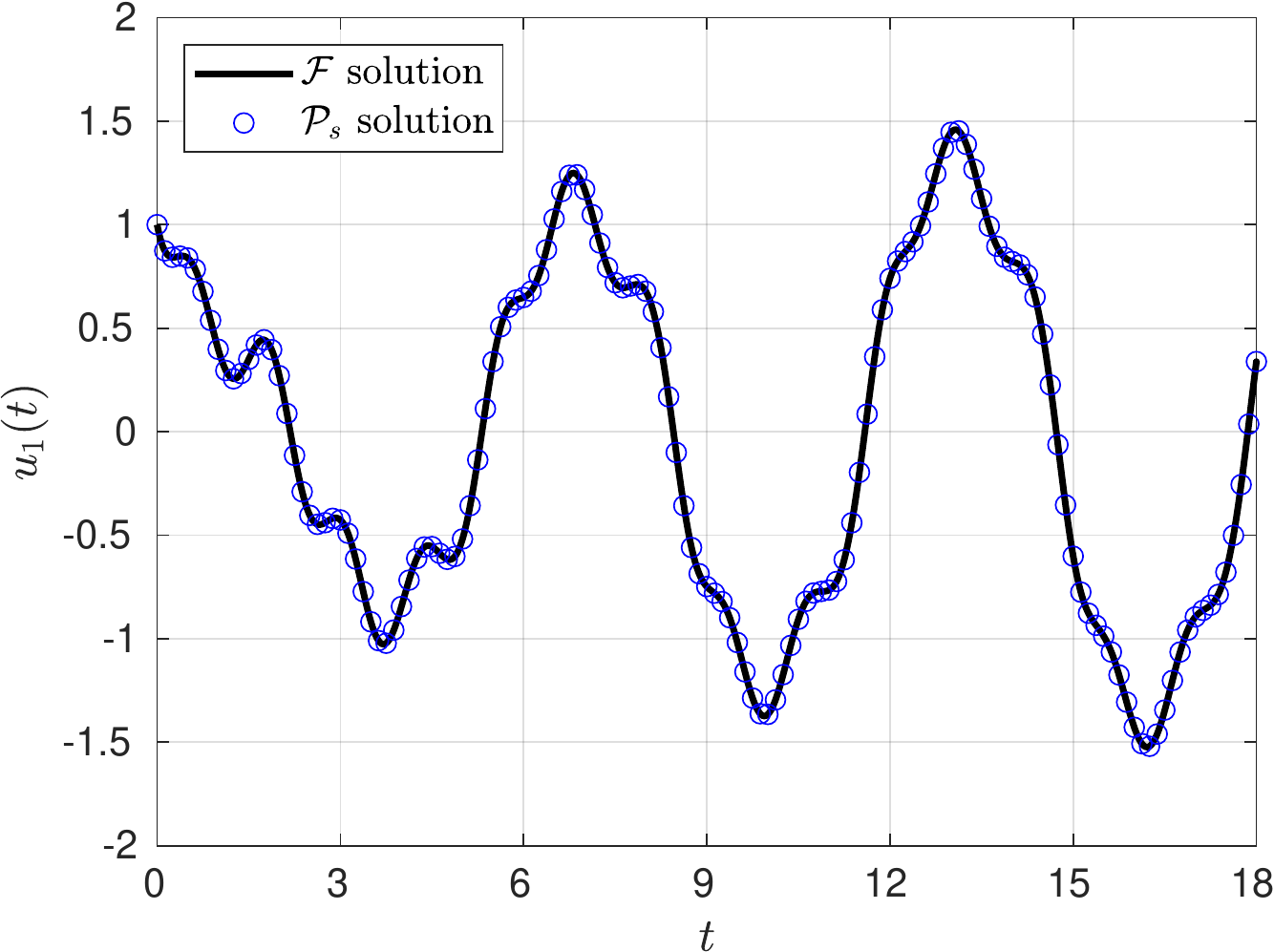} \label{fig:chartier_sol_a}}
    \subfloat[]{\includegraphics[width=0.49\textwidth]{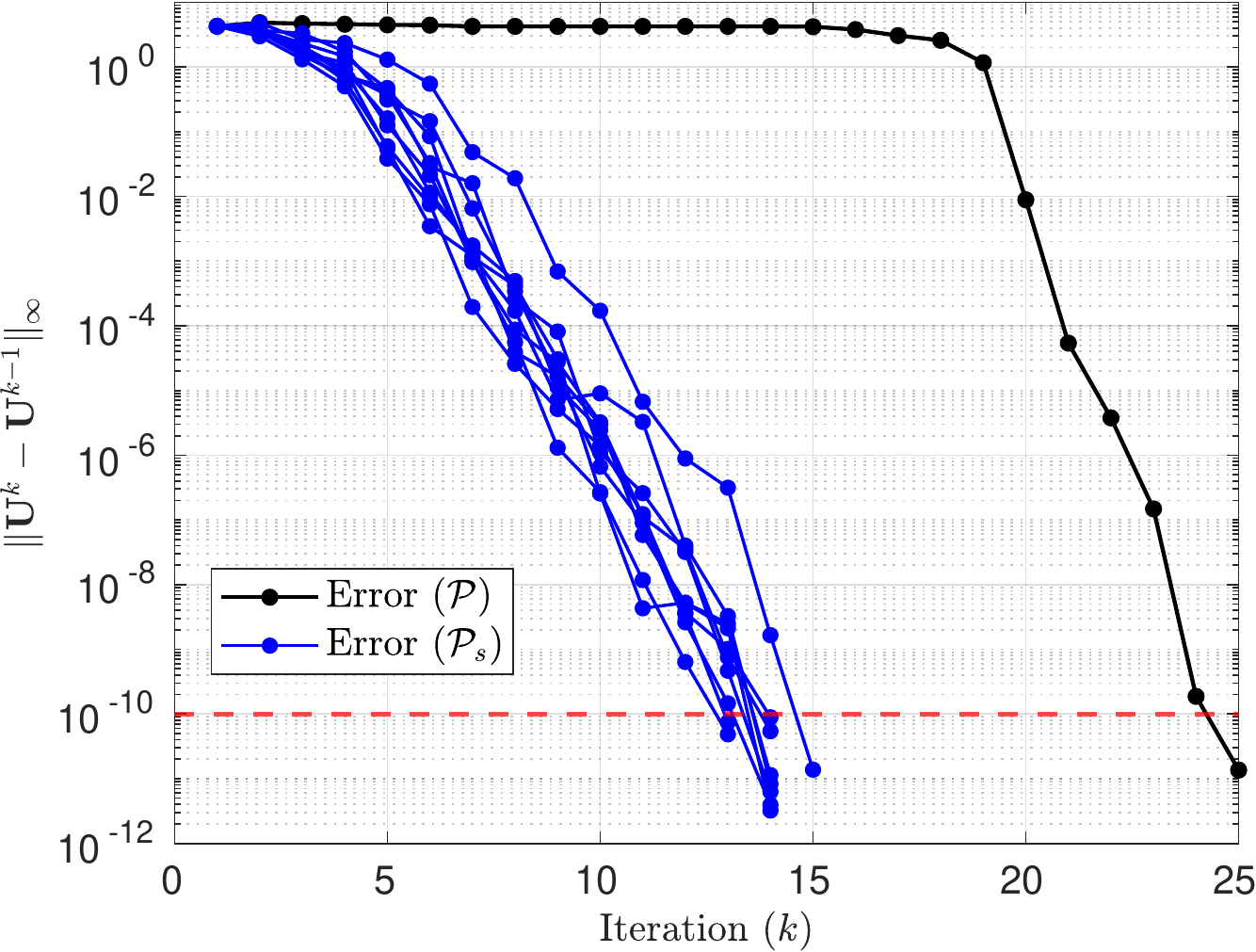} \label{fig:chartier_sol_b}}
    \caption{(a) Numerical solution of \eqref{eq:dissODE} over $[0,18]$ using $\mathcal{F}$ serially and a single realisation of $\mathcal{P}_s$. Note that only a subset of the fine times steps of the $\mathcal{P}_s$ solution are shown for clarity. (b) Errors at successive iterations of $\mathcal{P}$ (black line) and ten independent realisations of $\mathcal{P}_s$ (blue lines). Horizontal dashed red line represents the stopping tolerance $\varepsilon = 10^{-10}$. Note that both panels use $\mathcal{P}_s$ with sampling rule 1 and $M=3$.}
    \label{fig:chartier_sol}
\end{figure}

First, we consider the nonlinear ODE 
\begin{equation} \label{eq:dissODE}
    \frac{d u_1}{dt} = \sin(u_1) \cos(u_1) -2u_1 + e^{-t/100} \sin(5t) + \ln (1+t) \cos(t),
\end{equation}
with initial value $u_1(0) = 1$ \cite{Chartier1993}. Discretise the time interval $t \in [0,100]$ using $N=40$ sub-intervals, coarse time steps $\delta T = 100/80$, and fine time steps $\delta t = 100/8000$. Numerical solutions to \eqref{eq:dissODE} are shown on the interval $[0,18]$ in \Cref{fig:chartier_sol_a}. Deterministically, $\mathcal{P}$ locates a solution in $k_d = 25$ iterations using error tolerance $\varepsilon = 10^{-10}$.
$\mathcal{P}_s$ converges in a varying number of iterations $k_s$, with $\mathbb{P}(k_s<k_d) = 1$ -- see \Cref{fig:chartier_sol_b} for the convergence of ten independent simulations using $M=3$. 
From this plot, we can see that by taking just three samples, $\mathcal{P}_s$ reduces the convergence rate by almost a factor of two -- from $25$ to approximately $14$ (on average).
\begin{figure}[h]
    \centering
    \subfloat[]{\includegraphics[width=0.49\textwidth]{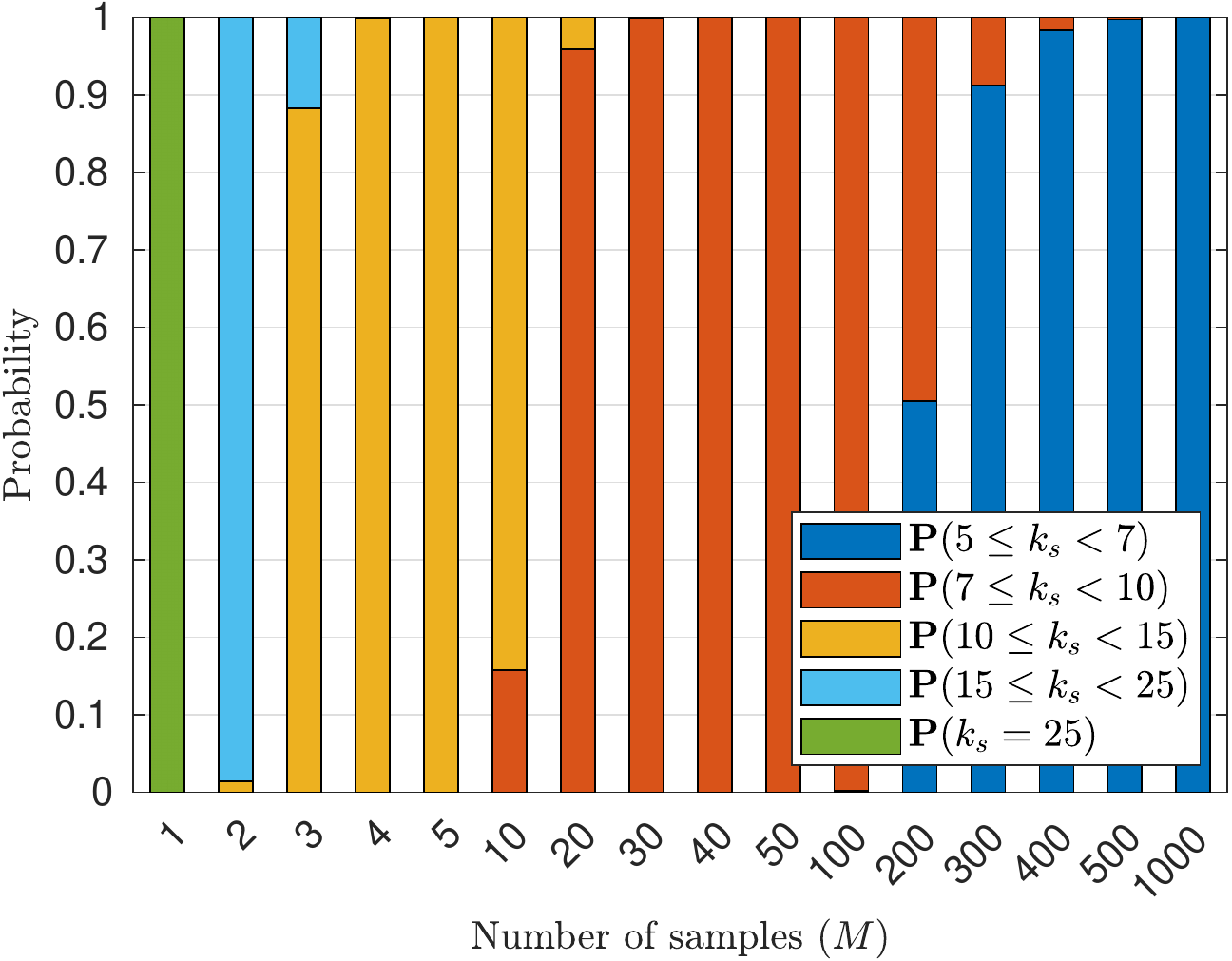} \label{fig:chartier_probs_a}}
    \subfloat[]{\includegraphics[width=0.49\textwidth]{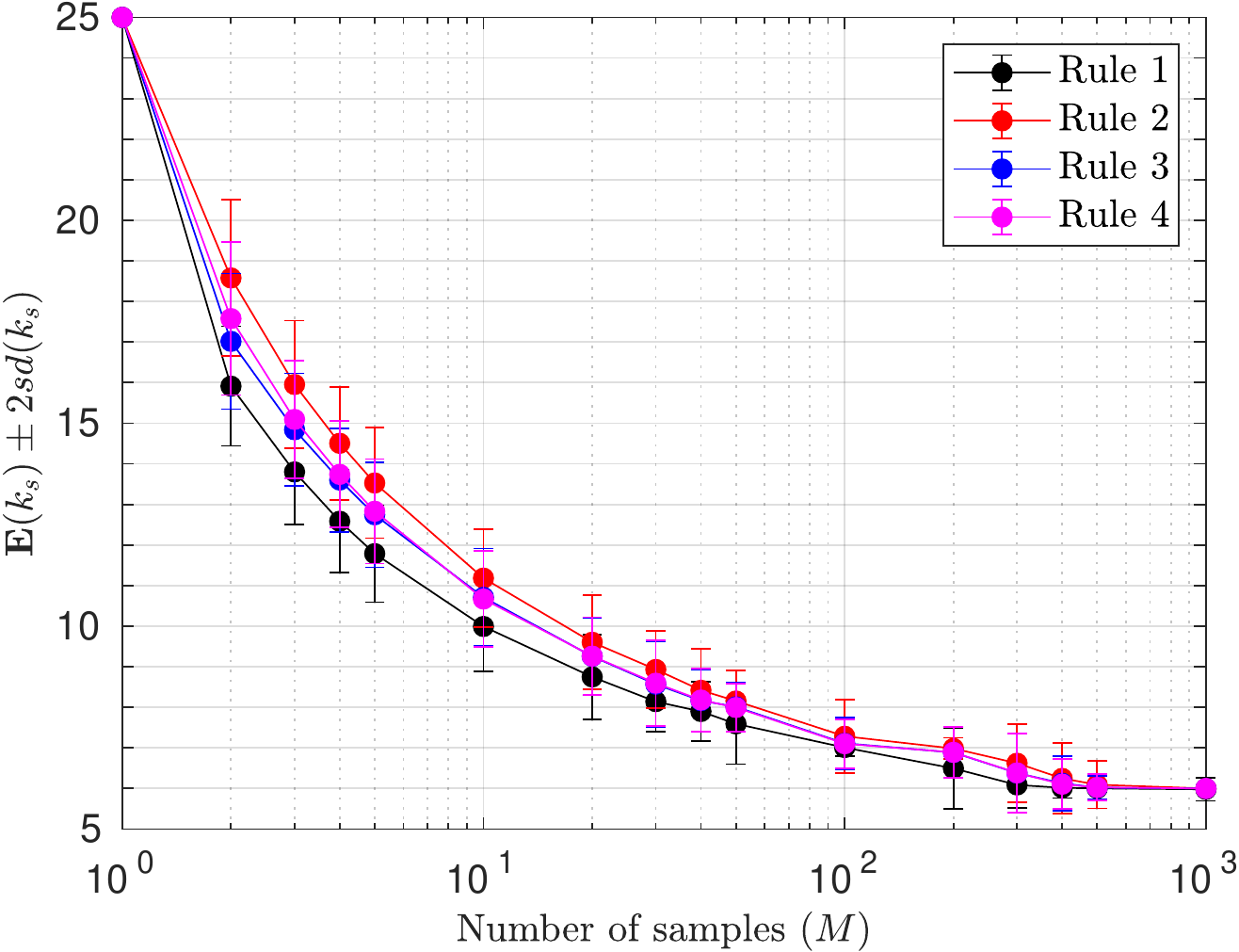} \label{fig:chartier_probs_b}}
    \caption{(a) Estimated discrete distributions of $k_s$ as a function of $M$ for sampling rule 1. (b) Estimated expectation of $k_s$ as a function of $M$, calculated using estimated distributions of $k_s$ for each sampling rule with error bars representing $\pm$ two standard deviations $sd(k_s)$. Distributions in both panels are estimated by simulating 2000 independent realisations of $\mathcal{P}_s$ for each $M$.}
    \label{fig:chartier_probs}
\end{figure}

When $M$ is increased above one, $\mathcal{P}_s$ begins generating stochastic solutions that converge in a varying number of iterations $k_s$. In order to accurately compare $k_d$ with the discrete random variable $k_s$, we run 2000 independent simulations of $\mathcal{P}_s$ to estimate the distribution of $k_s$ for a given $M$. Upon estimating these distributions, it was found that $\mathbb{P}(k_s < 25) = 1$ for each of the four sampling rules (for all $M>1$), meaning that by doubling the number of processors we can beat parareal with probability one. The estimated distributions of $k_s$, using sampling rule 1 (the other rules perform similarly), as a function of $M$ are given in \Cref{fig:chartier_probs_a}. The stacked bars represent the estimated discrete probability of a simulation converging in a given number of iterations. The results show $\mathcal{P}_s$ converging in just five iterations in the best case -- demonstrating $\mathcal{P}_s$ has the potential to yield significant parallel speedup, given sufficiently many samples are drawn. \Cref{fig:chartier_probs_b} emphasises the power of the stochastic method, showing that the estimated expected value $\mathbb{E}(k_s)$ decreases as $M$ increases, with the estimated standard deviation $sd(k_s) = \sqrt{\sum_{k=1}^N (k-\mathbb{E}(k_s))^2\mathbb{P}(k_s=k)}$ decreasing too. The improved performance of $\mathcal{P}_s$ as $M$ increases reflects what was discussed in \Cref{subsec:SP_convergence}. 
In particular, we ran a single realisation of $\mathcal{P}_s$ with $M=10^6$, observing that $\mathcal{P}_s$ converged in four iterations (result not shown), confirming that $k_s$ continues to decrease for increasing $M$.
By looking at \Cref{fig:chartier_probs_b}, we see that sampling rule 1 yields the lowest expected values of $k_s$ for small values of $M$, with all sampling rules performing similarly for large $M$.

\begin{figure}[t!]
    \centering
    \includegraphics[width=0.99\textwidth]{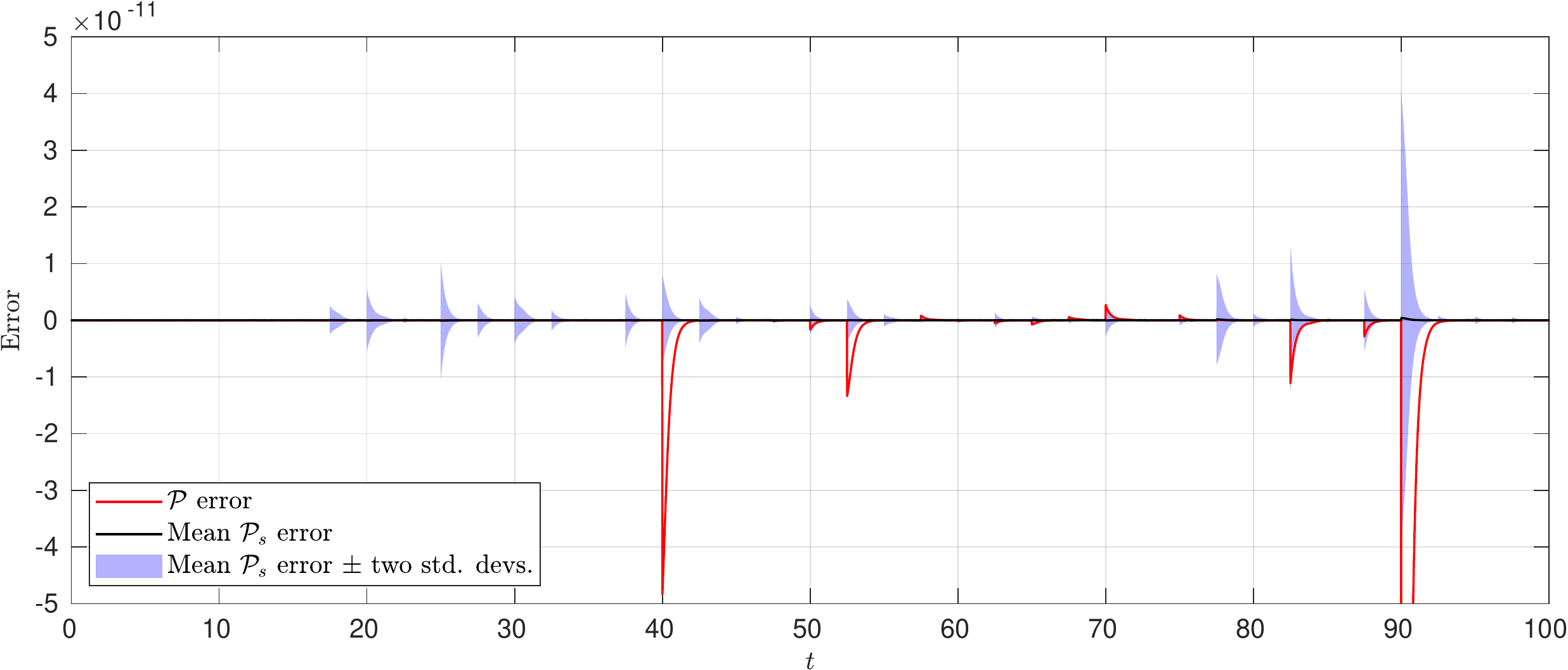}
    \caption{Errors of $\mathcal{P}$ (red) and mean $\mathcal{P}_s$ (black) solutions against the serial $\mathcal{F}$ solution over time. The mean error is obtained by running 2000 independent realisations of $\mathcal{P}_s$ with sampling rule 2 and $M=4$ -- the confidence interval representing the mean $\pm$ two standard deviations is shown in light blue.}
    \label{fig:chartier_errors}
\end{figure}

To verify the accuracy of the stochastic solutions, we plot the difference between the mean of 2000 independent realisations of $\mathcal{P}_s$ and the serially calculated $\mathcal{F}$ solution in \Cref{fig:chartier_errors}. Also shown is the confidence interval given by two standard deviations of the stochastic solutions (which is at most $\mathcal{O}(10^{-11})$) and the error generated by $\mathcal{P}$. Accuracy is maintained with respect to the fine solution across the time interval, even more so than the $\mathcal{P}$ solution. See \Cref{appen:1Dbern} for numerical results of $\mathcal{P}_s$ applied to a \textit{stiff} scalar nonlinear ODE. In this case, the stiffness of the equation demands a higher value of $M$ to improve $k_s$ -- something we also observe for the Brusselator in \Cref{subsec:2Dbruss}.

%%%%%%%%%%%%%%%%%%%%%%%%%%%%%%%%%%%%%%%%%%%%%%%%%%%%%%%%%%%%%%%%%
\subsection{The Brusselator system} \label{subsec:2Dbruss}
Next, consider the Brusselator system
\begin{subequations} \label{eq:brussODE}
\begin{align} 
        \frac{d u_1}{dt} &= A + u_1^2 u_2 - (B+1)u_1, \\
        \frac{d u_2}{dt} &=  Bu_1 - u_1^2 u_2,
\end{align}    
\end{subequations}
a pair of stiff nonlinear ODEs that model an auto-catalytic chemical reaction \cite{Lefever1971}. Using parameters $(A,B) = (1,3)$, trajectories of the system exhibit oscillatory behaviour in phase space, approaching a limit cycle (as $t \rightarrow \infty$) that contains the unstable fixed point $(1,3)^{\text{T}}$. Now that $d>1$, we use bivariate distributions to sample the initial values -- meaning we can compare the effects of including or excluding the correlations between variables. System \eqref{eq:brussODE} is solved using initial values $\mathbf{u}(0) = (1,3.07)^{\text{T}}$ over time interval $t \in [0,15.3]$ with $N=25$, $\delta T = 15.3/25$, and $\delta t = 15.3/2500$ \cite{Trefethen2017}. The numerical solution to \eqref{eq:brussODE} in phase space and convergence of the successive errors are reported in \Cref{appen:2Dbruss}. With these parameters and a tolerance of $\varepsilon = 10^{-6}$, $\mathcal{P}$ takes $k_d = 7$ iterations to stop and return a numerical solution. 

\begin{figure}[b!]
    \centering
    \subfloat[]{\includegraphics[width=0.48\textwidth]{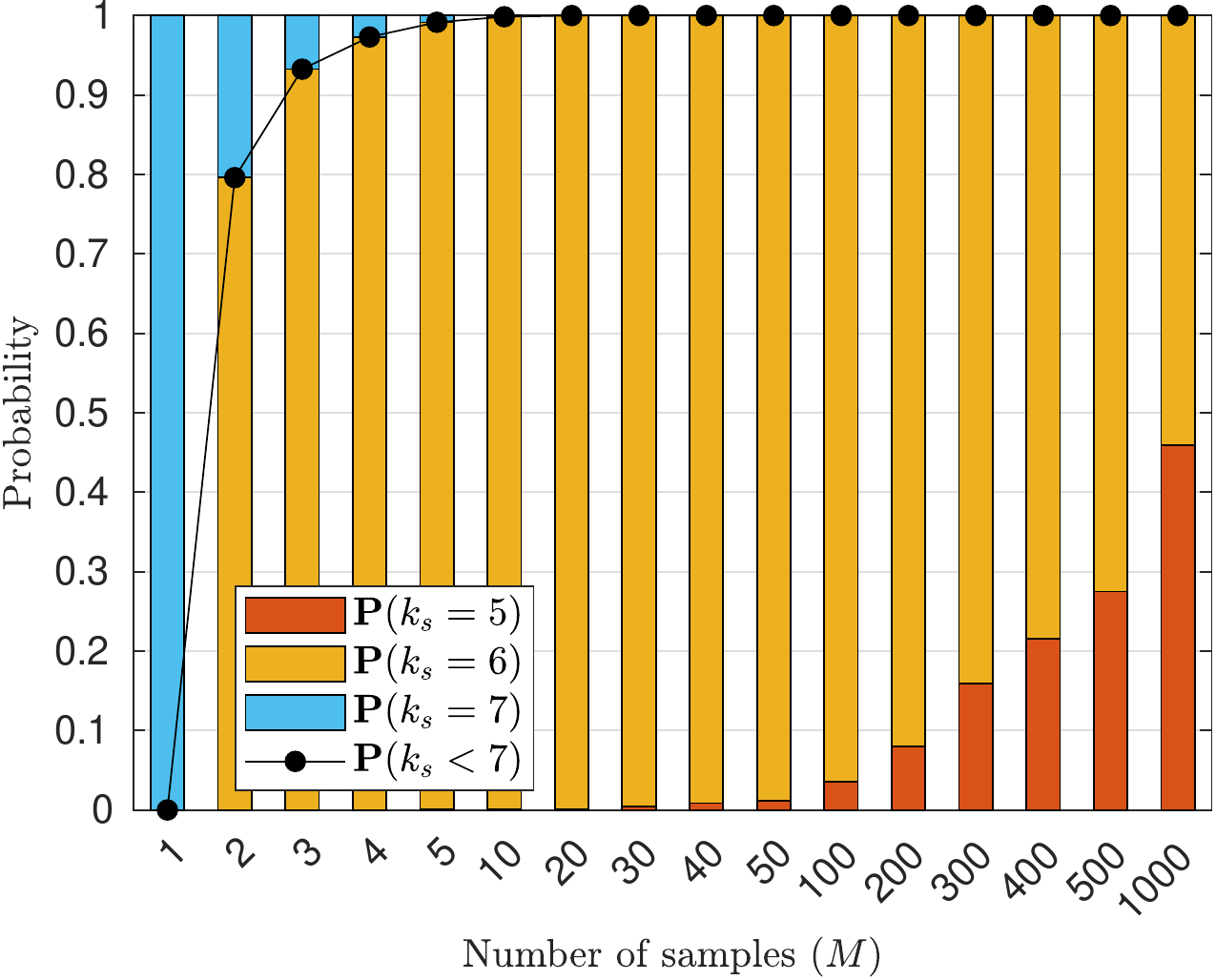} \label{fig:bruss_probs_a}}
    \subfloat[]{\includegraphics[width=0.49\textwidth]{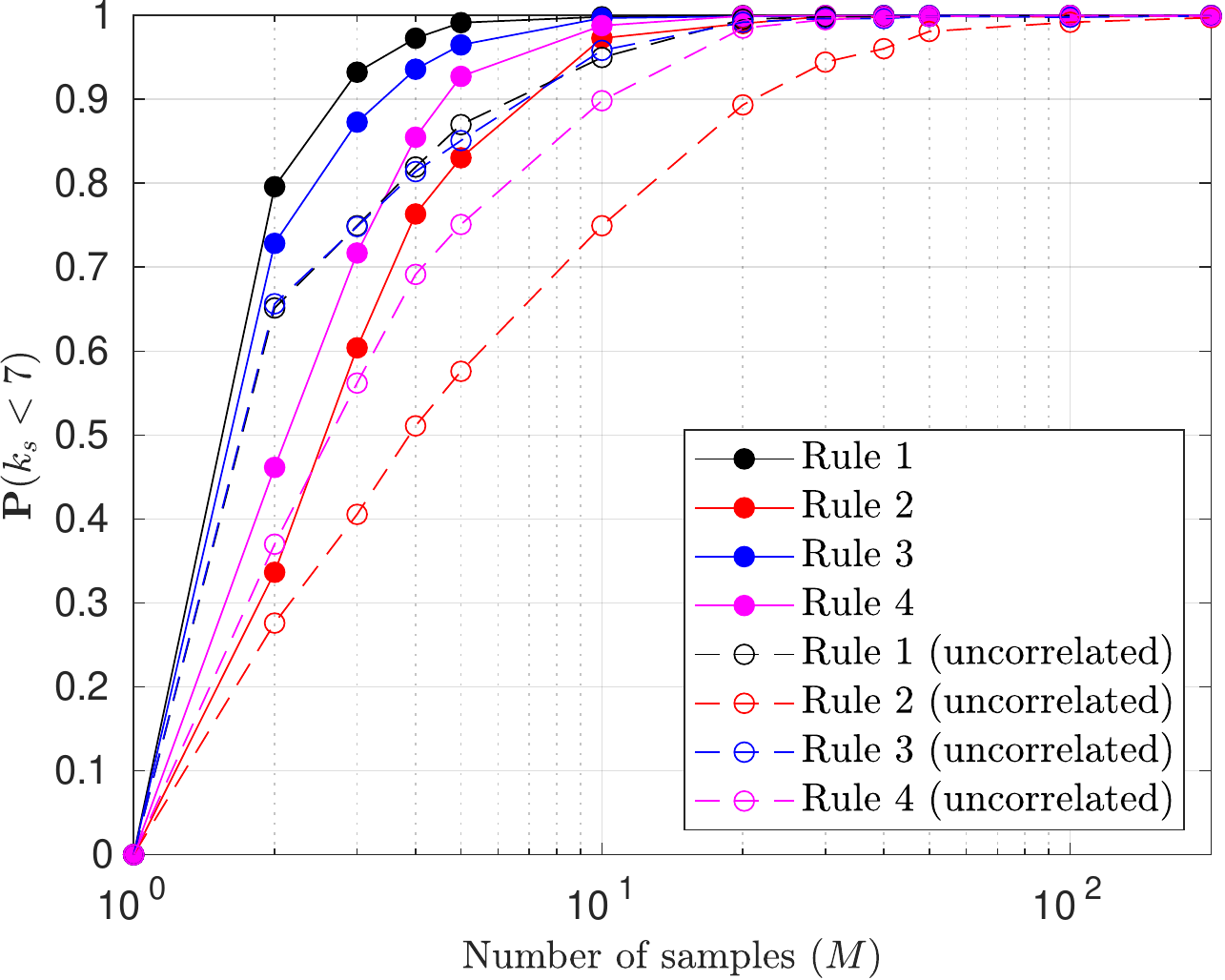} \label{fig:bruss_probs_b}}
    \caption{(a) Estimated discrete probabilities of $k_s$ as a function of $M$ for sampling rule 1. (b) Estimated probability that the convergence rate $k_s$ is smaller than $k_d=7$ as a function of $M$ for the sampling rules with (solid lines) and without (dashed lines) correlations. Distributions were estimated by simulating 2000 independent realisations of $\mathcal{P}_s$ for each $M$.}
    \label{fig:bruss_probs}
\end{figure}

\begin{figure}[t!]
    \centering
    \subfloat[]{\includegraphics[width=0.49\textwidth]{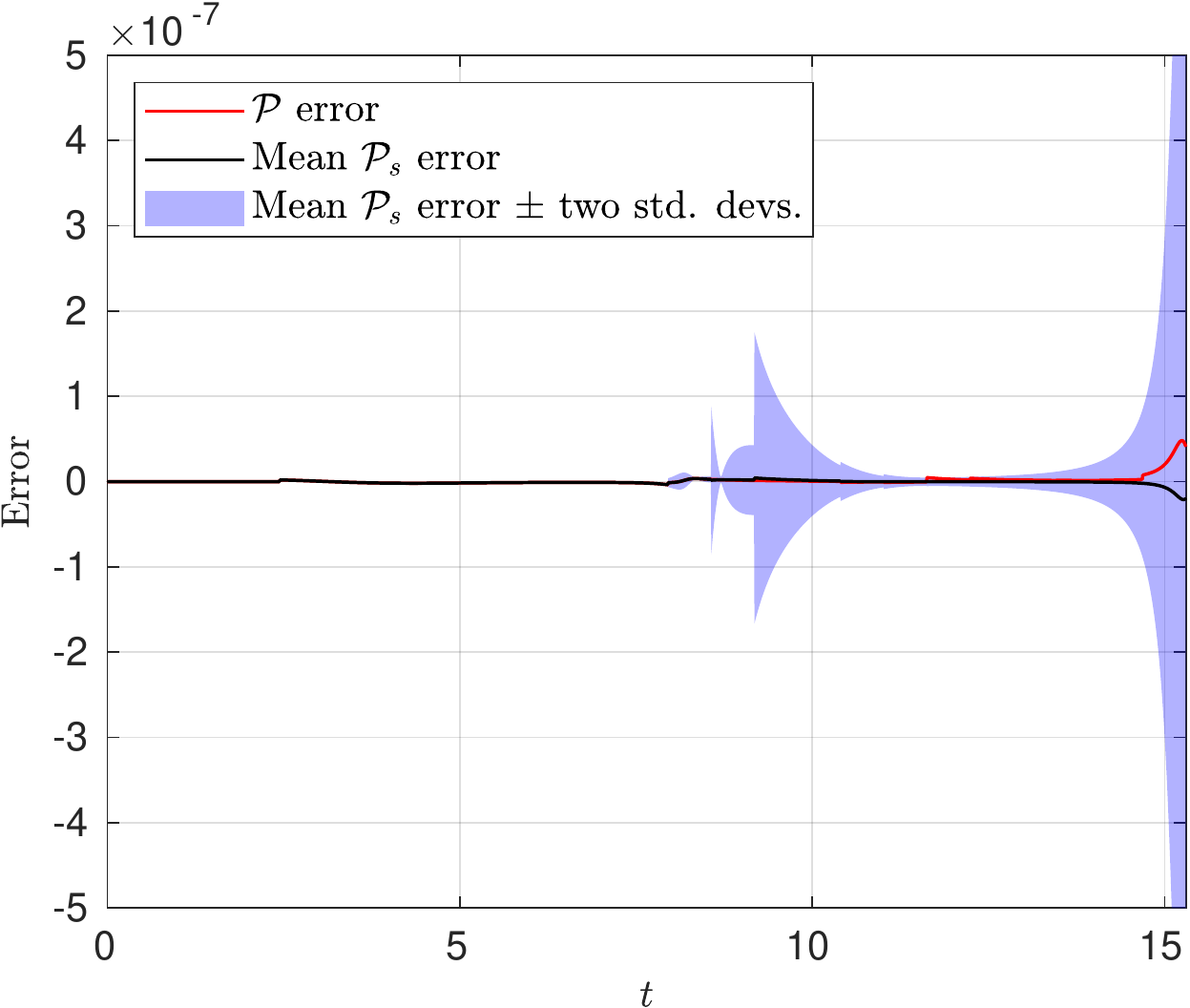} \label{fig:bruss_errors_a}}
    \subfloat[]{\includegraphics[width=0.49\textwidth]{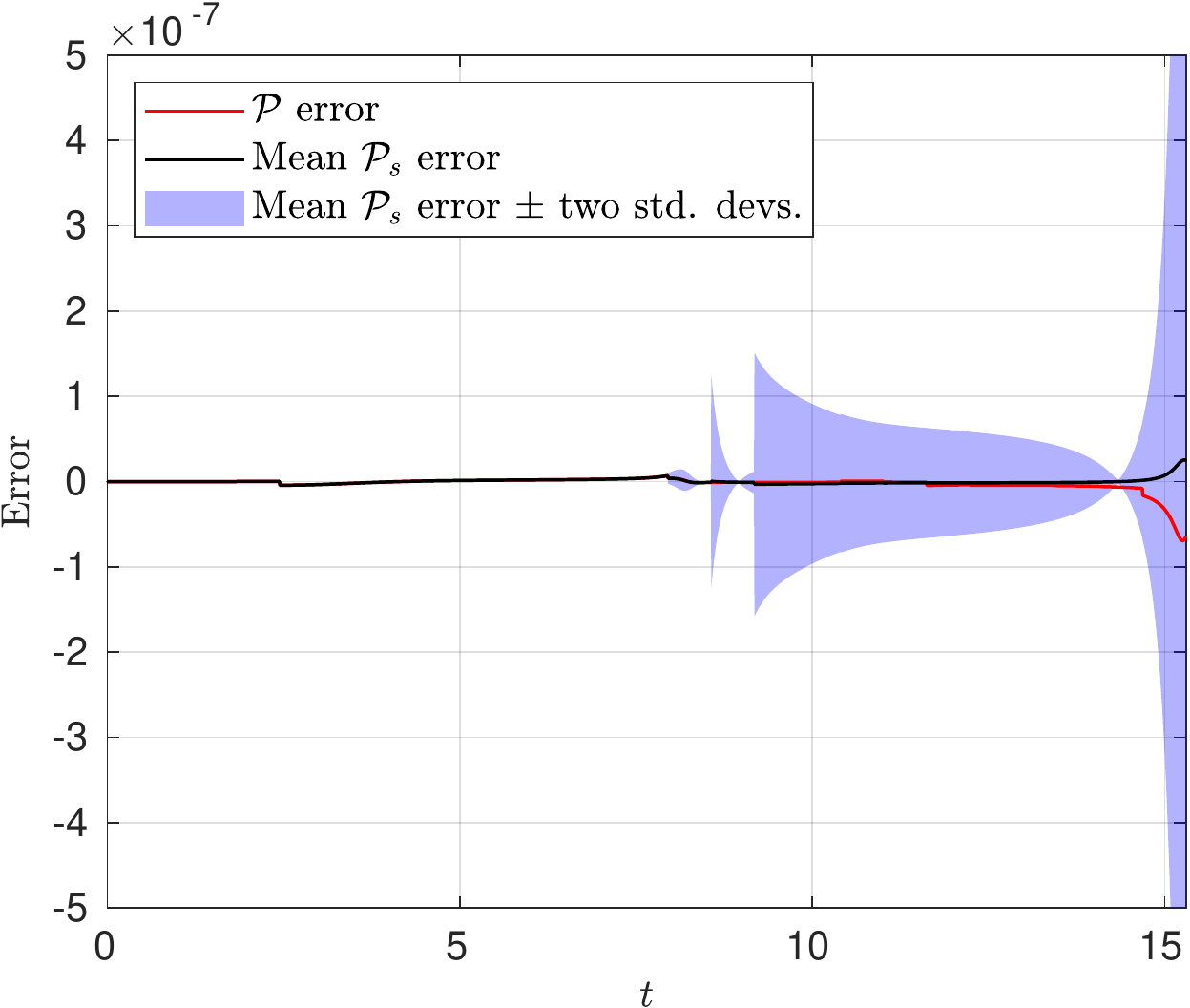} \label{fig:bruss_errors_b}}        
    \caption{Errors of $\mathcal{P}$ (red) and mean $\mathcal{P}_s$ (black) solutions against the $\mathcal{F}$ solution. The mean error is obtained by running 2000 independent realisations of $\mathcal{P}_s$ with sampling rule 4 with $M=200$ -- the confidence interval representing the mean $\pm$ two standard deviations is shown in light blue. Panel (a) displays errors for the $u_1$ component of the solution whilst (b) displays the $u_2$ component.}
    \label{fig:bruss_errors}
\end{figure}

The estimated distributions of $k_s$ for sampling rule 1 are given in \Cref{fig:bruss_probs_a}. Even though $\mathcal{P}$ takes just $k_d=7$ iterations to stop, we observe that $\mathcal{P}_s$ can still reach the desired tolerance in 5 or 6 iterations -- albeit requiring larger values of $M$. We believe this is due to the stiffness of the system and poor accuracy of the explicit $\mathcal{G}$ solver -- results presented for the stiff ODE in \Cref{appen:1Dbern} appear to confirm this. Using adaptive time-stepping methods could be a way to reduce the value of $M$ needed to converge in fewer $k_s$ \cite{Maday2020}. The solid lines in \Cref{fig:bruss_probs_b} show that, using sampling rules 1 or 3, $\mathcal{P}_s$ only requires $M \approx 10$ to beat parareal almost certainly, i.e.\ to guarantee that $\mathbb{P}(k_s<7) \rightarrow 1$. Sampling rules 1 and 3 outperform 2 and 4 in this particular system. Note, however, the stark decrease in performance if instead uncorrelated samples are generated within $\mathcal{P}_s$ (dashed lines). This demonstrates the importance of accounting for the dependence between variables in nonlinear systems such as \eqref{eq:brussODE}. Observe again, in \Cref{fig:bruss_errors}, how the mean $\mathcal{P}_s$ solutions attain equivalent, or better, accuracy than the $\mathcal{P}$ solutions, with standard deviations at most $\mathcal{O}(10^{-6})$.

Further results of $\mathcal{P}_s$ applied to \eqref{eq:brussODE} and an additional two-dimensional nonlinear system are presented in \Cref{appen:2Dbruss} and \Cref{appen:2Dsquare} respectively. Observe again that for the less stiff system in \Cref{appen:2Dsquare}, we require less sampling to improve $k_s$ compared to the Brusselator -- further highlighting the demand for higher $M$ in stiff systems to improve the convergence rate.

%%%%%%%%%%%%%%%%%%%%%%%%%%%%%%%%%%%%%%%%%%%%%%%%%%%%%%%%%%%%%%%%%
\subsection{The Lorenz system} \label{subsec:3Dlorenz}
Finally, we consider the Lorenz system
\begin{subequations} \label{eq:lorenzODE}
\begin{align} 
        \frac{d u_1}{dt} &= \gamma_1 (u_2 - u_1), \\
        \frac{d u_2}{dt} &= \gamma_2 u_1 - u_1 u_3 - u_2, \\
        \frac{d u_3}{dt} &= u_1 u_2 - \gamma_3 u_3,
\end{align}    
\end{subequations}
a simplified model for weather prediction \cite{Lorenz1963}. With the parameters $(\gamma_1,\gamma_2,\gamma_3) = (10,28,8/3)$, \eqref{eq:lorenzODE} exhibits chaotic behaviour where trajectories with initial values close to one another diverge exponentially. This will test the robustness of $\mathcal{P}_s$, as small numerical differences between initial values will mean that errors can grow rapidly as time progresses. We solve \eqref{eq:lorenzODE} using initial values $\mathbf{u}(0) = (-15,-15,20)^{\text{T}}$ over the interval $[0,18]$, discretised using $N=50$ sub-intervals and time steps $\delta T = 18/250$ and $\delta t = 18/18750$. With a tolerance of $\varepsilon = 10^{-8}$, $\mathcal{P}$ takes $k_d=20$ iterations to converge.
\begin{figure}[t!]
    \centering
    \subfloat[]{\includegraphics[width=0.48\textwidth]{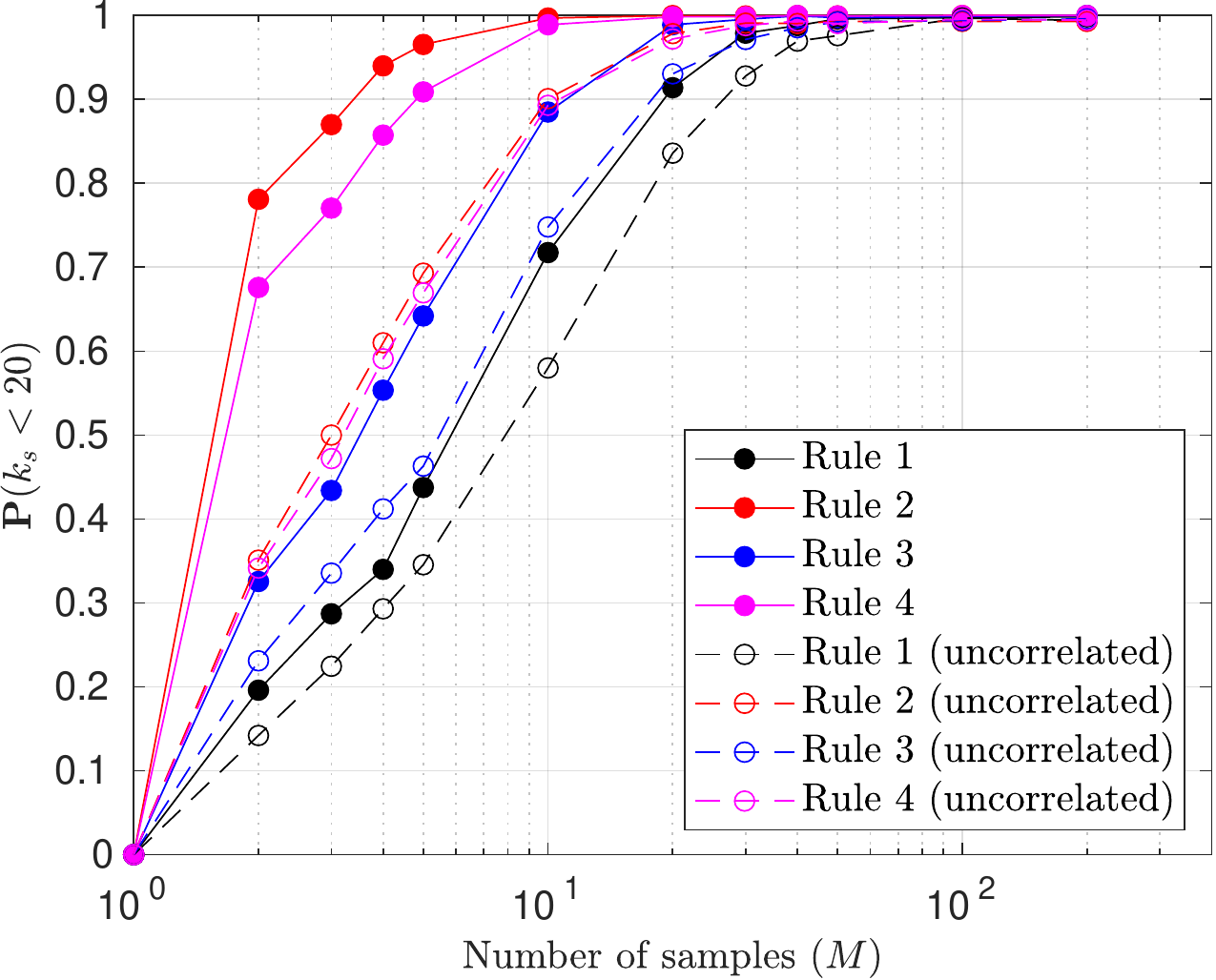} \label{fig:lorenz_probs_a}}
    \subfloat[]{\includegraphics[width=0.49\textwidth]{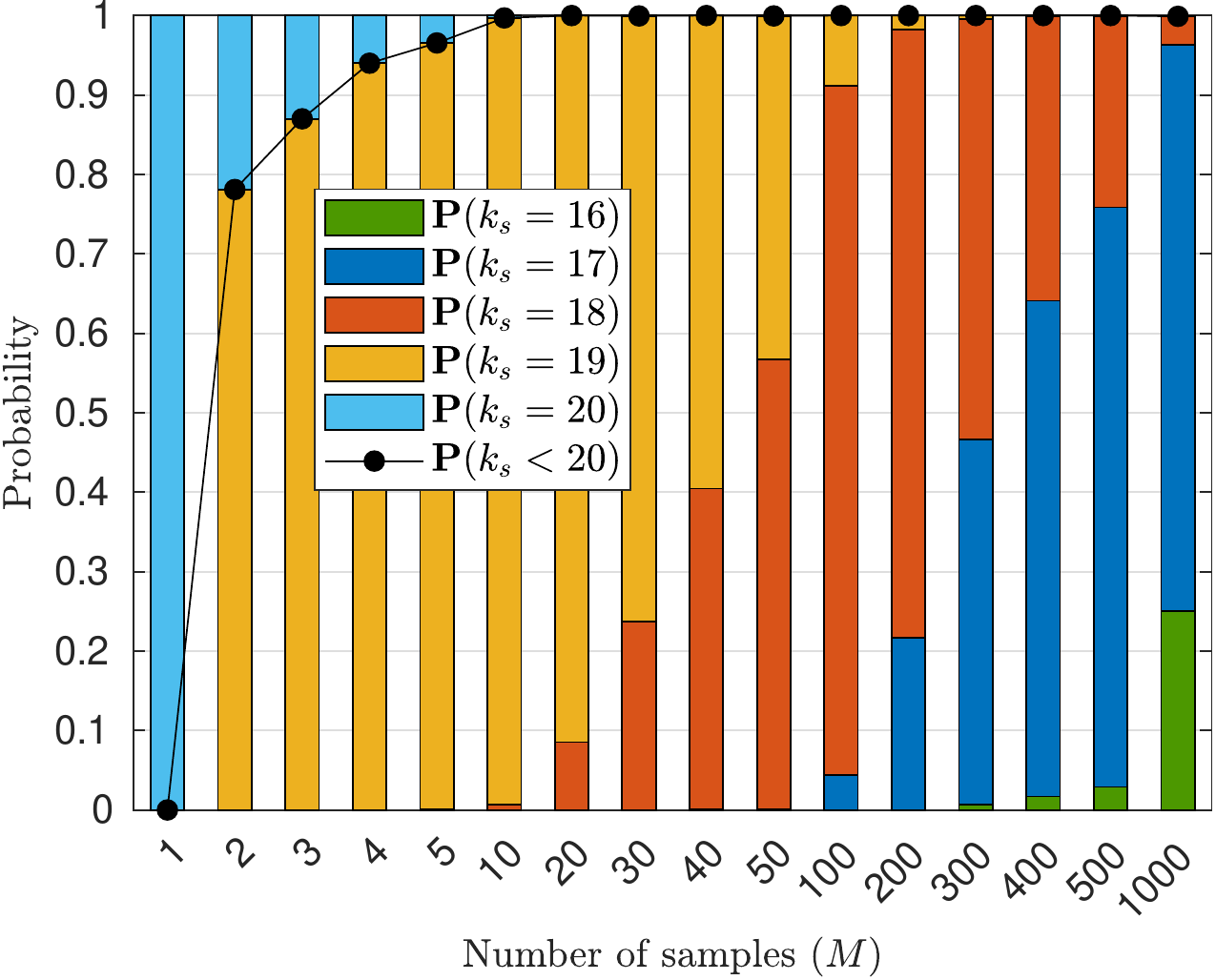} \label{fig:lorenz_probs_b}}
    \caption{(a) Estimated probability that the convergence rate $k_s$ is smaller than $k_d=20$ as a function of $M$ for the sampling rules with (solid lines) and without (dashed lines) correlations. Distributions were estimated by simulating 2000 independent realisations of $\mathcal{P}_s$ for each $M$. (b) Estimated discrete probabilities of $k_s$ as a function of $M$ for sampling rule 2. }
    \label{fig:lorenz_probs}
\end{figure}

Running $\mathcal{P}_s$ to compare the performance of the sampling rules, we see again in \Cref{fig:lorenz_probs_a} that taking correlated samples is much more efficient than not and that only $M \approx 10$ samples are required to beat parareal with probability one. For the chaotic trajectories generated by \eqref{eq:lorenzODE}, sampling close to the predictor-corrector, rules 2 and 4, yields superior performance compared to rules 1 and 3 for small values of $M$. \Cref{fig:lorenz_probs_b} displays estimated distributions for varying $M$ using sampling rule 2 -- yielding a best convergence rate $k_s = 16$ for approximately 25\% of runs with $M=1000$. These results demonstrate the robustness of $\mathcal{P}_s$ and that the sampling and propagation process is not impeded by the exponential divergence of trajectories. Additional results of $\mathcal{P}_s$ applied to \eqref{eq:lorenzODE} are given in \Cref{appen:3Dlorenz}.

%% file: sec5.tex
% !TeX root = ms.tex

\section{Conclusions} \label{sec:conclusions}

%what did we do and what are the benefits of stochastic parareal?
In this paper, we have extended the parareal algorithm using probabilistic methods to develop a stochastic parareal algorithm for solving systems of ODEs in a time-parallel manner. 
Instead of passing deterministically calculated initial values into parareal's predictor-corrector, stochastic parareal selects more accurate values from a randomly sampled set, at each temporal sub-interval, to converge in fewer iterations. 
In \Cref{sec:results}, we compared performance against the deterministic parareal algorithm on several low-dimensional ODE systems of increasing complexity by calculating the estimated convergence rate distributions (upon multiple independent realisation of stochastic parareal) with increasing numbers of random samples $M$.
By taking just $M\approx10$ (correlated) samples, the estimated probability of converging sooner than parareal approached one in all test cases.
Similarly, we observed numerical convergence toward the fine (exact) solution with accuracy of similar order to parareal and, in the spirit of probabilistic numerics \cite{hennig2015,oates2019}, obtained a measure of uncertainty over the ODE solution upon multiple realisations of the algorithm.

The probability that stochastic parareal converges faster than standard parareal depends on a number of factors: the complexity of the problem being solved, the number of time sub-intervals ($N$), the accuracy of the coarse integrator ($\mathcal{G}$), the number of random samples ($M$), and the type of sampling rule in use. 
Sampling rules 1 and 3 (sampling close to the fine solutions) outperformed rules 2 and 4 (sampling close to the predictor-corrector solutions) for the ODE systems in \Cref{subsec:1Dnonlinear}, \Cref{subsec:2Dbruss}, and SM1. 
The reverse was true, however, for the systems in \Cref{subsec:3Dlorenz} and SM3, making it difficult to determine an optimal rule for a general ODE system. 
To overcome having to choose a particular sampling rule, one could linearly combine different rules or even sample from multiple rules simultaneously. 
We would suggest sampling from probability distributions with infinite support, i.e. the Gaussians (rules 1 and 2), so that samples can be taken anywhere in $\mathbb{R}^d$ with non-zero probability.
Having finite support may have created difficulty for the uniform marginal $t$-copulas (rules 3 and 4) because samples could only be taken in a finite hyperrectangle in $\mathbb{R}^d$ -- problematic if the exact solution were to lay outside of this space.

%To tackle this, one could consider new sampling distributions using Gaussian or $t$-copulas with underlying marginal distributions of infinite support, i.e.\ defined on $\mathbb{R}$.

When solving stiff ODEs (see \Cref{subsec:2Dbruss} and SM1), results indicated that stochastic parareal demanded increasingly high sampling to converge sooner than parareal than for non-stiff systems.
For example, we observe that when taking $M=100$ samples in the non-stiff scalar ODE in \Cref{fig:chartier_probs}, the expected convergence rate decreases from 25 to 7 whereas for the stiff scalar ODE in SM1, the rate only drops from 8 to 6. 
A similar observation can be made in the two-dimensional test cases in \Cref{subsec:2Dbruss} (stiff) and SM3 (non-stiff). These results exemplify the role that system complexity, e.g. stiffness or chaos, plays in the performance of both algorithms.
In SM1, stochastic parareal was also shown to perform more efficiently for problems that deterministic parareal itself struggles with, i.e. cases in which the accuracy of the coarse integrator $\mathcal{G}$ is poor. 
In SM3 it was also observed that, for low sample numbers ($M = 2$), stochastic parareal actually converged in one more iteration than parareal in less than $2.5\%$ of cases. This suggests there may be minimum number of samples required to beat parareal in some situations -- something to be investigated with further experimentation. 

%future work
In summary, we have demonstrated that probabilistic methods and additional processors can, for low-dimensional ODEs, be used to increase the parallel scalability of existing time-parallel algorithms. 
Next we need to investigate possible improvements and generalisations. 
For example, determining whether the algorithm scales for larger systems of equations -- essential if it is to be used for solving PDE problems.
Moreover, whether we can design adaptive sampling rules that do not need to be specified \textit{a priori} to simulation. 
Finally, we would aim to make use of the whole ensemble of fine propagated trajectories rather than using only one -- avoiding  the waste of valuable information about the solution at the coarse and fine resolutions. 
An approach in this direction has recently been proposed \cite{pentland2022}, making use of ideas from the field of probabilistic numerics to adopt a more Bayesian approach to this problem.

%% file: sec6.tex
% !TeX root = ms.tex

In the following sections we illustrate some additional results using stochastic parareal.

\section{Scalar Bernoulli equation} \label{appen:1Dbern}
\begin{figure}[b!]
    \centering
    \subfloat[]{\includegraphics[width=0.49\textwidth]{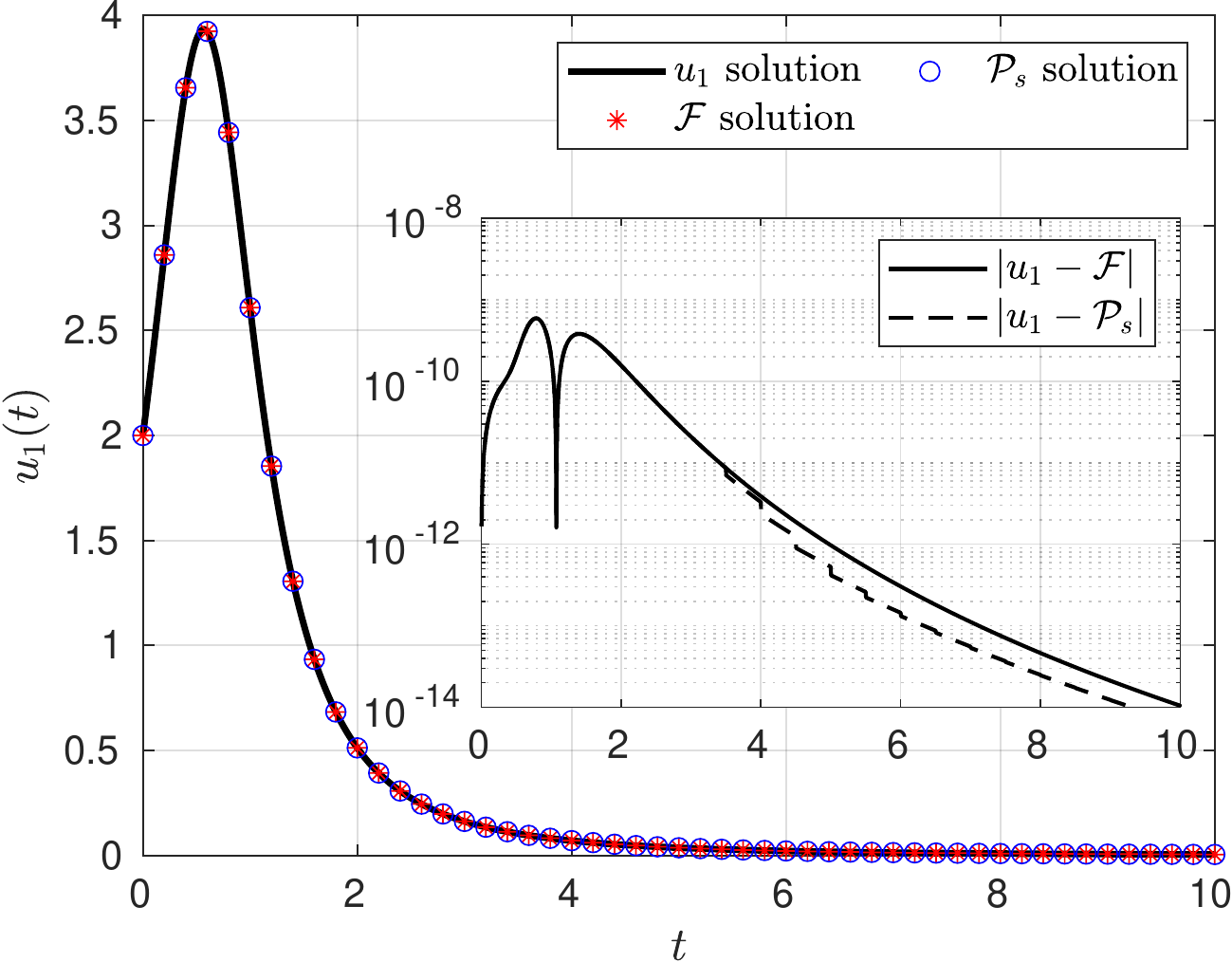} \label{fig:bern_sol_a}}
    \subfloat[]{\includegraphics[width=0.49\textwidth]{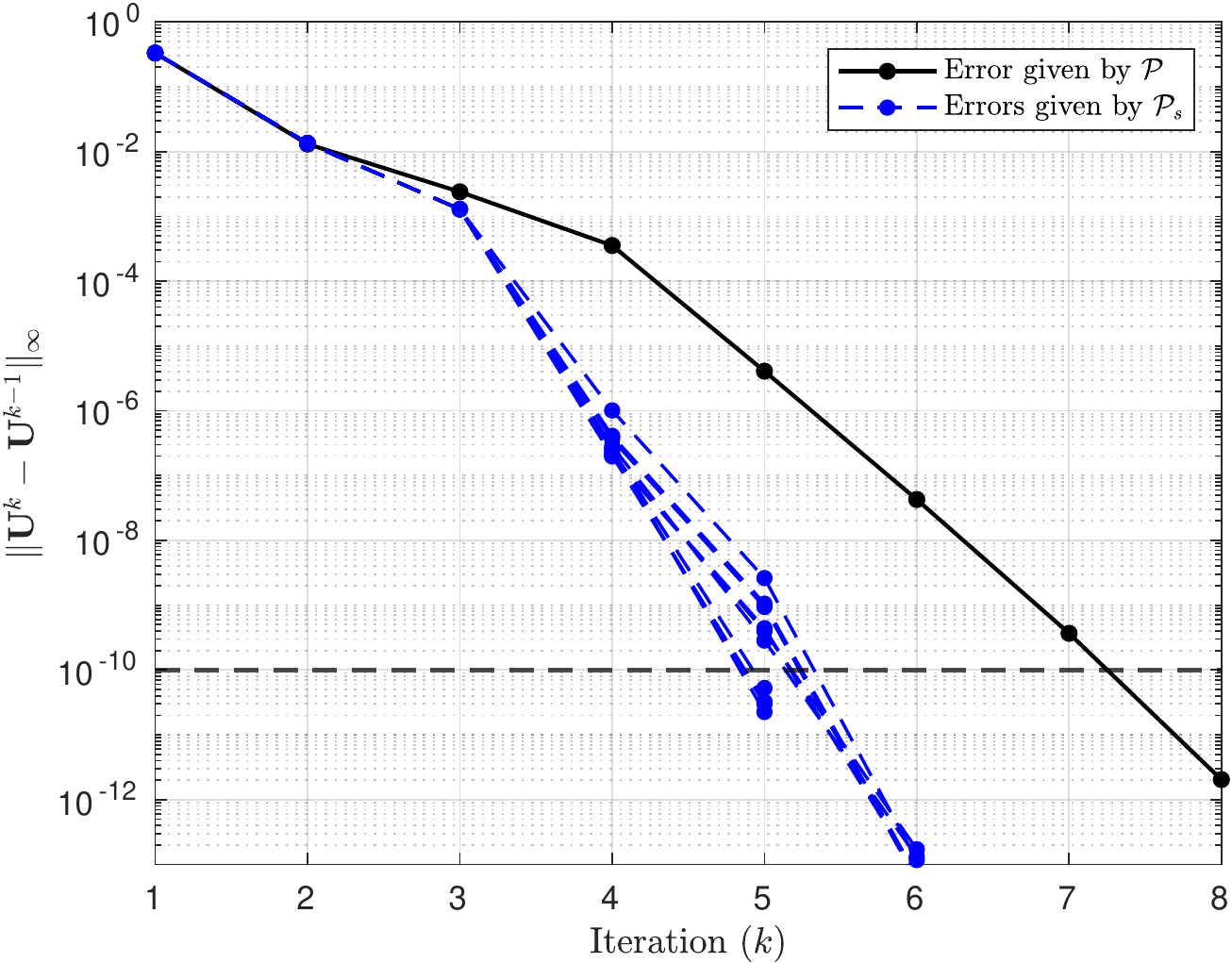} \label{fig:bern_sol_b}}
    \caption{(a) Analytical solution $u_1$ of \eqref{eq:bernODE} plotted against the serial $\mathcal{F}$ solution and a single realisation of $\mathcal{P}_s$. Note that only a subset of the fine times steps of the numerical solutions are shown for clarity. Inset: displays the numerical errors of $\mathcal{F}$ and $\mathcal{P}_s$ compared to $u_1$. (b) Errors at successive iterations of $\mathcal{P}$ (black lines) and ten independent realisations of $\mathcal{P}_s$ (blue lines). The horizontal dashed black line represents the tolerance $\varepsilon = 10^{-10}$. Both panels use $\mathcal{P}_s$ with sampling rule 3 and (a) $M=30$ and (b) $M=1000$.}
    \label{fig:bern_sol}
\end{figure}

Consider the nonlinear non-autonomous Bernoulli equation
\begin{equation} \label{eq:bernODE}
    \frac{d u_1}{dt} = \frac{2}{1+t} u_1 - t^2 u_1^2,
\end{equation}   
with initial value $u_1(0) = 2$ on $t \in [0,10]$. We discretise using $N=20$ sub-intervals and time steps $\delta T = 10/20$ and $\delta t = 10/2000$. Equation \eqref{eq:bernODE} permits the analytical solution $u_1(t) = (1+t)^2 / (t^5 /5 + t^4 /2 + t^3 / 3 + 1/2)$, tending to zero as $t \rightarrow \infty$. Observe the spacing between equidistant time steps of the true $\mathcal{F}$ solution and the $\mathcal{P}_s$ solution to \eqref{eq:bernODE} in \cref{fig:bern_sol_a}, highlighting the stiffness of the solution at early times. Given the stopping tolerance $\varepsilon = 10^{-10}$, observe in \cref{fig:bern_sol_b} how $\mathcal{P}$ converges in $k_d = 8$ iterations deterministically whilst $\mathcal{P}_s$ converges in $k_s \in \{5,6 \}$ iterations for each of the ten independent realisations shown with $M=1000$.
\begin{figure}[htbp]
    \centering
    \subfloat[]{\includegraphics[width=0.49\textwidth]{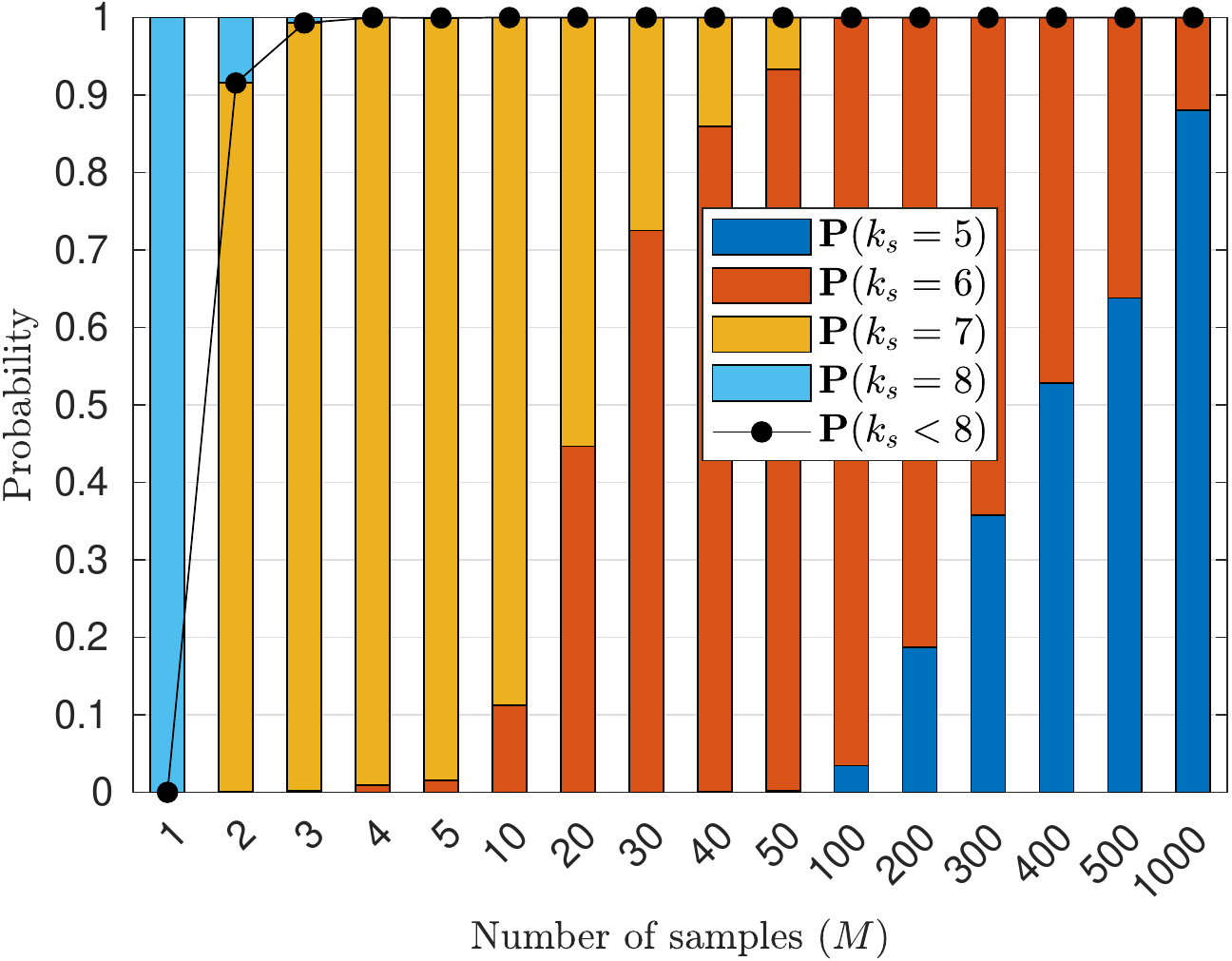} \label{fig:bern_probs_a}}
    \subfloat[]{\includegraphics[width=0.49\textwidth,height=0.39\textwidth]{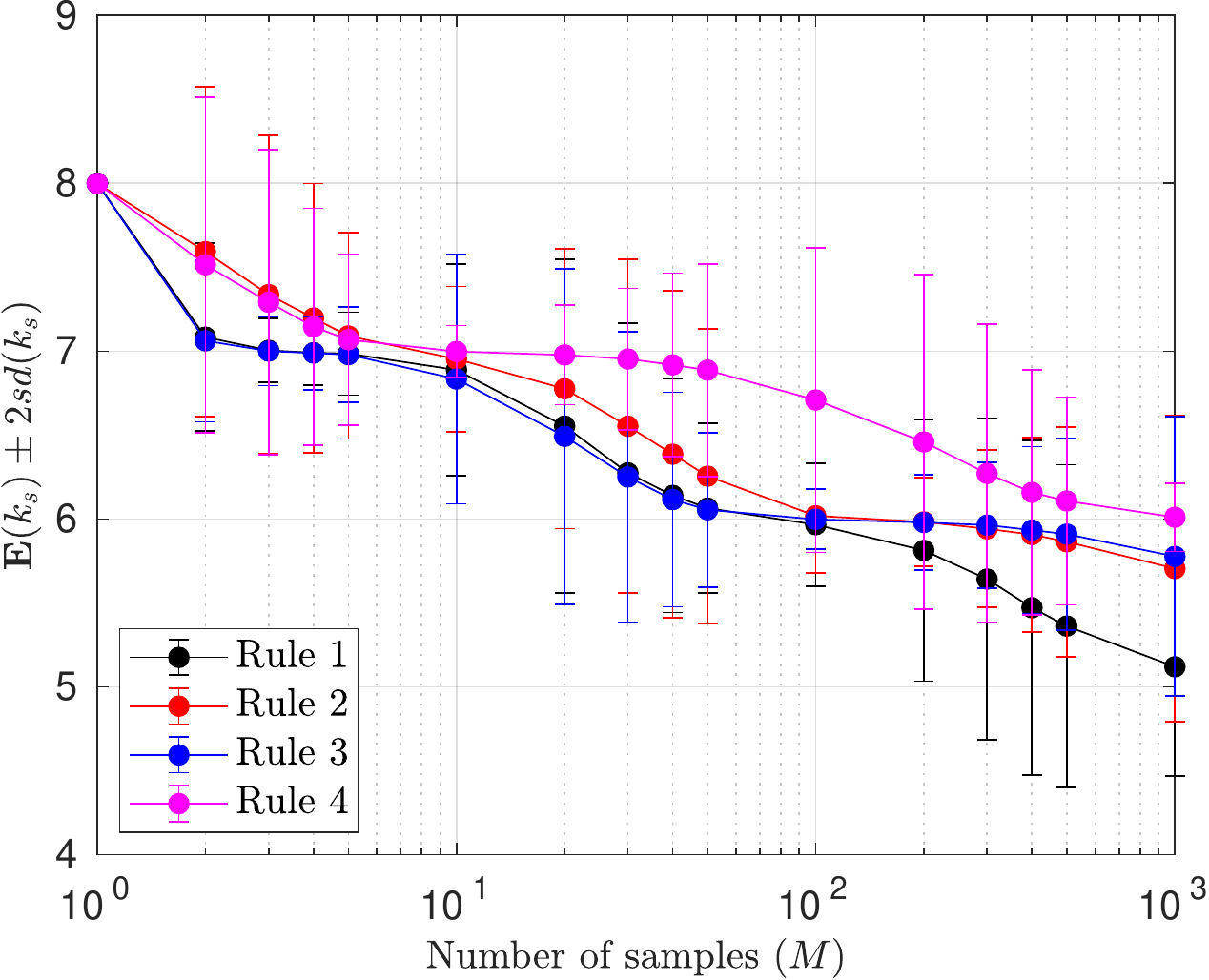} \label{fig:bern_probs_b}}    
    \caption{(a) Estimated discrete distributions of $k_s$ as a function of $M$ for sampling rule 1. (b) Estimated expectation of $k_s$ as a function of $M$, calculated using estimated distributions of $k_s$ for each sampling rule with error bars representing $\pm$ two standard deviations $sd(k_s)$. Distributions were estimated by simulating 2000 independent realisations of $\mathcal{P}_s$ for each $M$.}
    \label{fig:bern_probs}
\end{figure}

\Cref{fig:bern_probs_a} shows the estimated distributions of $k_s$ using sampling rule 1. As expected, if $M=1$, convergence is deterministic (i.e. $\mathcal{P}_s = \mathcal{P}$) and hence $\mathbb{P}(k_s = 8) = 1$.  As $M$ increases however, $\mathbb{P}(k_s = 8)$ decreases rapidly to zero whilst $\mathbb{P}(k_s < 8)$ increases from zero to one with just $M \approx 5$ samples. This demonstrates that $\mathcal{P}_s$ requires very few samples to begin converging in fewer iterations than $\mathcal{P}$, and that $k_s$ assumes low values with increasing probabilities for increasing $M$. The stiffness of \eqref{eq:bernODE} appears to demand much larger values of $M$ to continually reduce $k_s$ when compared to the non-stiff scalar example in subsection 4.1.

\begin{figure}[htbp]
    \centering
    \includegraphics[width=0.48\textwidth]{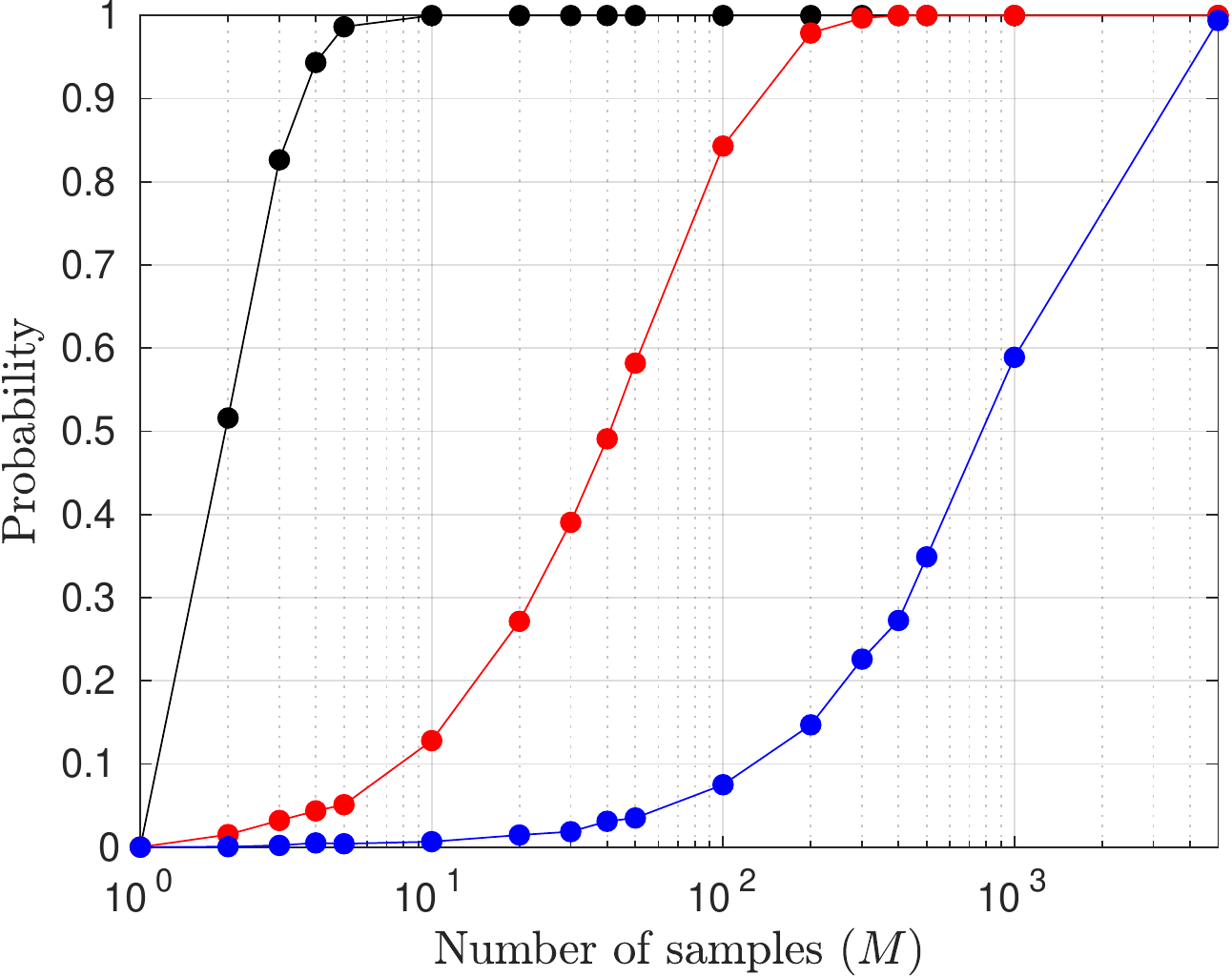}
    \caption{Estimated probabilities that the convergence rate $k_s$ is smaller than $k_d$ as a function of $M$ using sampling rule 3. Each curve shows $\mathbb{P}(k_s < 8)$ (black), $\mathbb{P}(k_s < 5)$ (red), and $\mathbb{P}(k_s < 4)$ (blue) using coarse steps $\delta T = 10/20, 10/40, 10/60$ respectively - noting that $\mathcal{P}$ converges in $k_d = 8, 5 , 4$ iterations for these coarse steps respectively. As before, 2000 independent realisations of $\mathcal{P}_s$ were run for each $M$.}
    \label{fig:bern_coarse_probs}
\end{figure}

In \cref{fig:bern_probs_b} we report the expected values of $k_s$ for each sampling rule. This shows that for low values of $M$, using either a Gaussian or $t$-copula sampling rule makes little difference to performance. More interestingly, rules 1 and 3, centred around the fine solutions, outperform rules 2 and 4, which are centred around the predictor-corrector. Further testing revealed that varying the fine time steps within $\mathcal{P}_s$ had little impact on these probabilities. On the contrary, \cref{fig:bern_coarse_probs} shows how increasing the number of coarse steps in the interval $[0,10]$ from $20$ to $60$ drastically decreases the probabilities of $\mathcal{P}_s$ converging sooner than $\mathcal{P}$. Increasing the number of coarse steps increases the accuracy of the $\mathcal{G}$ solver, hence $\mathcal{P}$ reaches the stopping tolerance in fewer iterations $k_d$. This result suggests that whilst $\mathcal{P}_s$ can still converge faster than $\mathcal{P}$ by using more samples, it works more efficiently for particular problems where $k_d$ is relatively large (i.e when $\mathcal{P}$ obtains low rates of parallel speed up). 

\begin{figure}[htbp]
    \centering
    \subfloat[]{\includegraphics[width=0.49\textwidth]{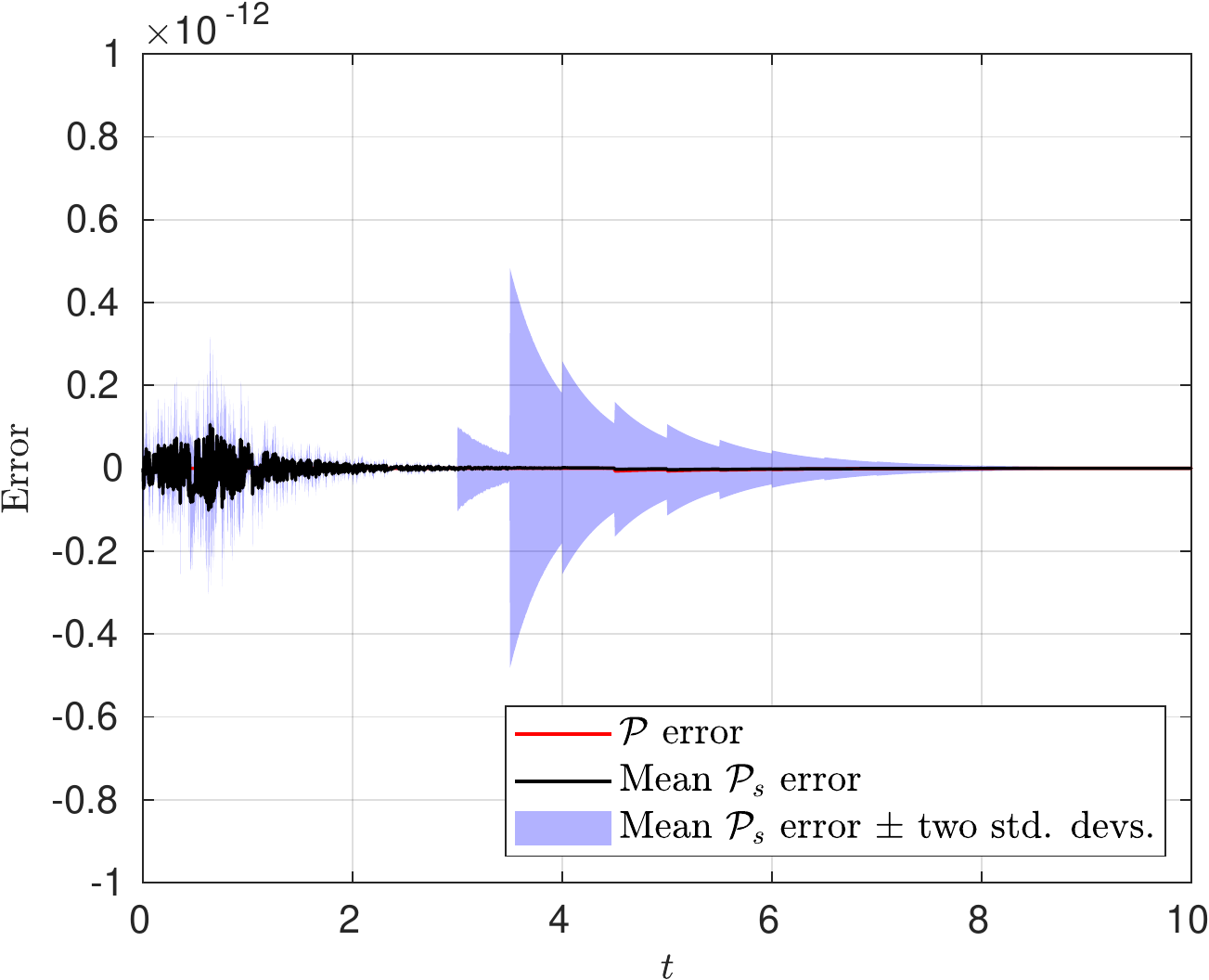} \label{fig:bern_errors_a}}
    \subfloat[]{\includegraphics[width=0.49\textwidth]{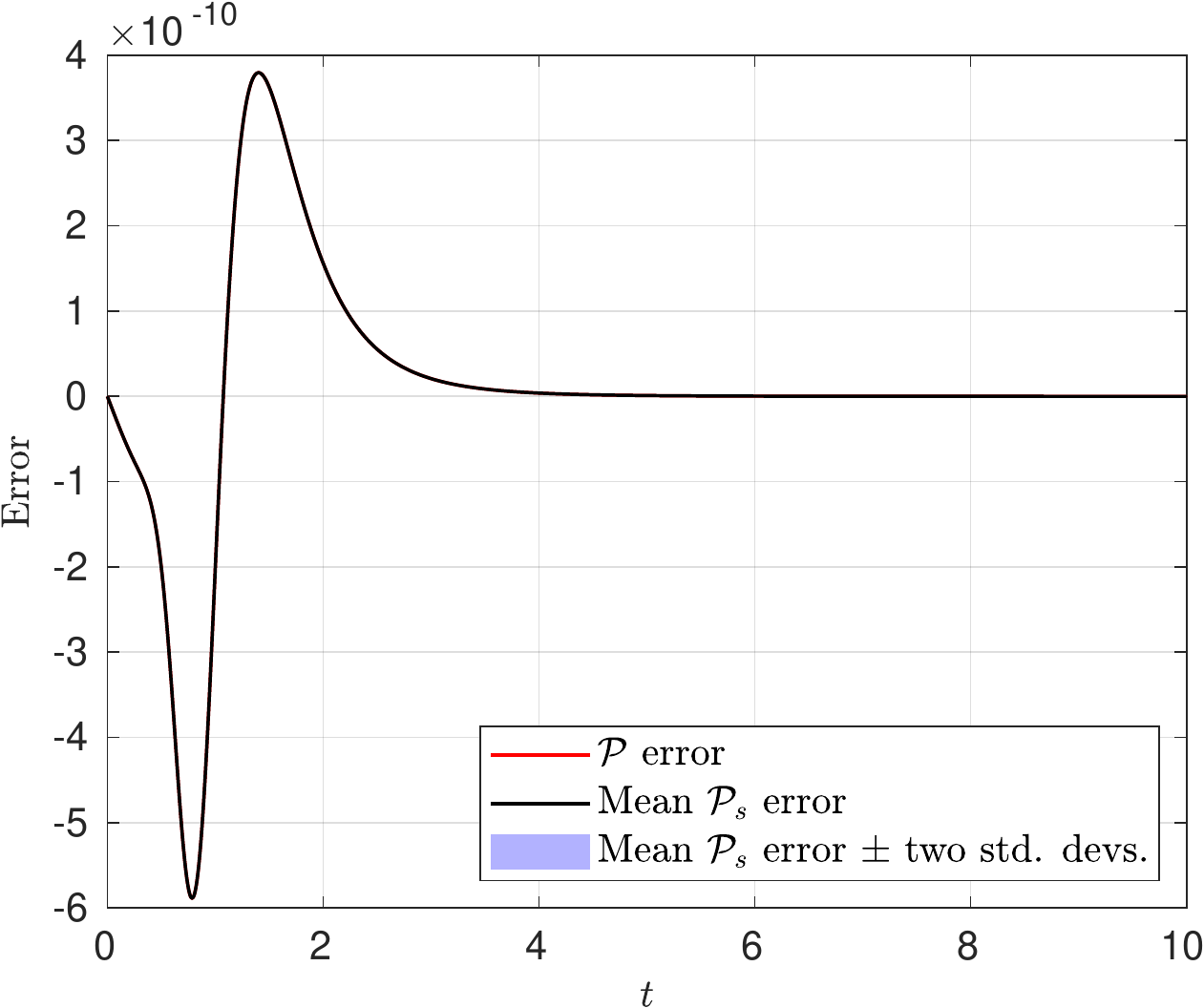} \label{fig:bern_errors_b}}        
    \caption{Errors of $\mathcal{P}$ (red) and mean $\mathcal{P}_s$ (black) solutions against (a) the $\mathcal{F}$ solution and (b) the analytical solution $u_1$. The mean error is obtained by running 2000 independent realisations of $\mathcal{P}_s$ with sampling rule 1 with $M=10$ -- the confidence interval representing the mean $\pm$ two standard deviations is shown in light blue.}
    \label{fig:bern_errors}
\end{figure}

Finally, we calculated errors between the mean $\mathcal{P}_s$ solution with respect to the $\mathcal{F}$ solution (\cref{fig:bern_errors_a}) and the analytical solution $u_1$ (\cref{fig:bern_errors_b}). In both cases, we observe how the mean solution attains comparable accuracy to $\mathcal{P}$ and both the fine and analytical solutions.

%%%%%%%%%%%%%%%%%%%%%%%%%%%%%%%%%%%%%%%%%%%%%%%%%%%%%%%%%%%%%%%%%
\newpage
\section{The Brusselator system} \label{appen:2Dbruss}

\begin{figure}[htbp]
    \centering
    \subfloat[]{\includegraphics[width=0.49\textwidth]{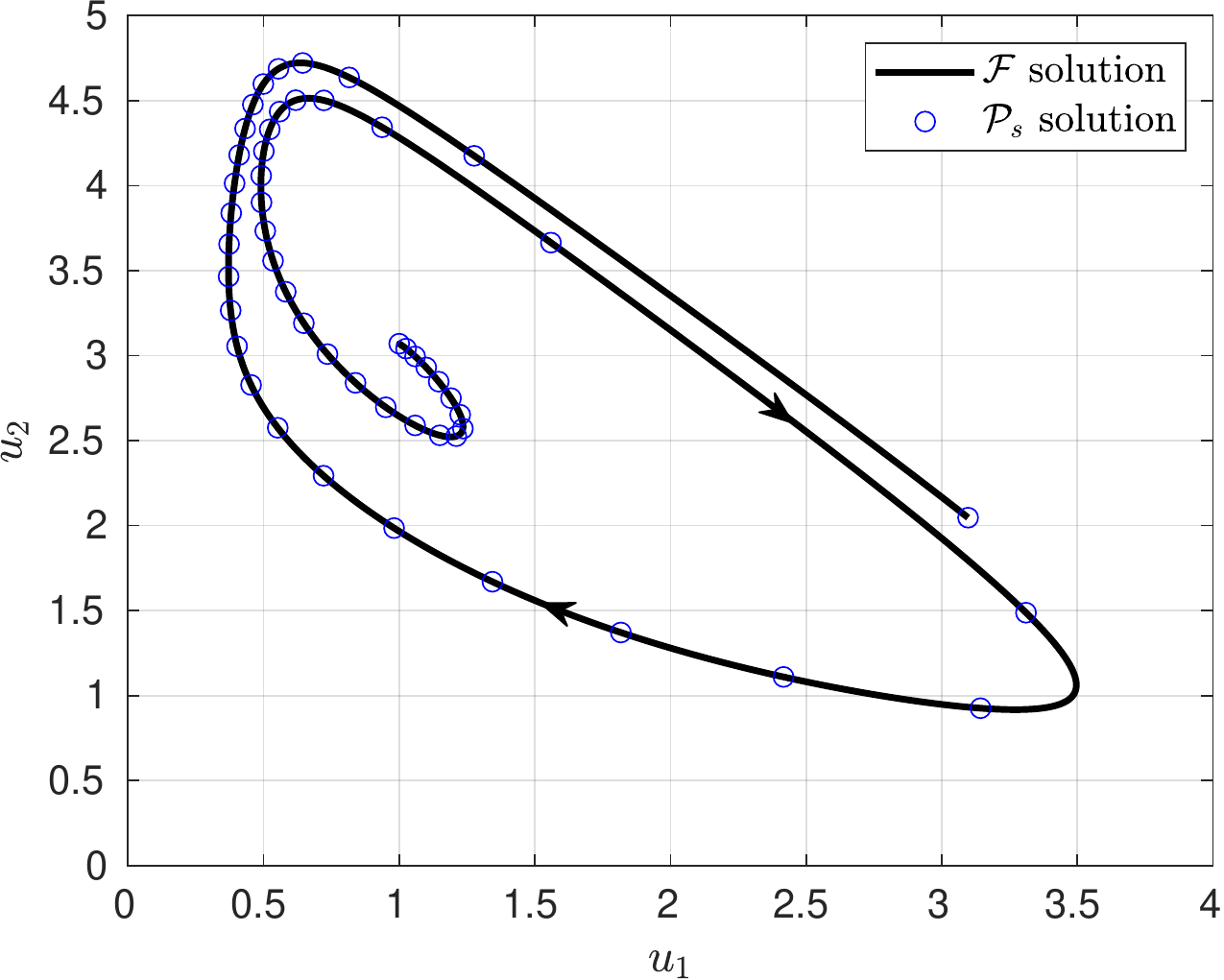} \label{fig:bruss_sol_a}}
    \subfloat[]{\includegraphics[width=0.49\textwidth,height=0.385\textwidth]{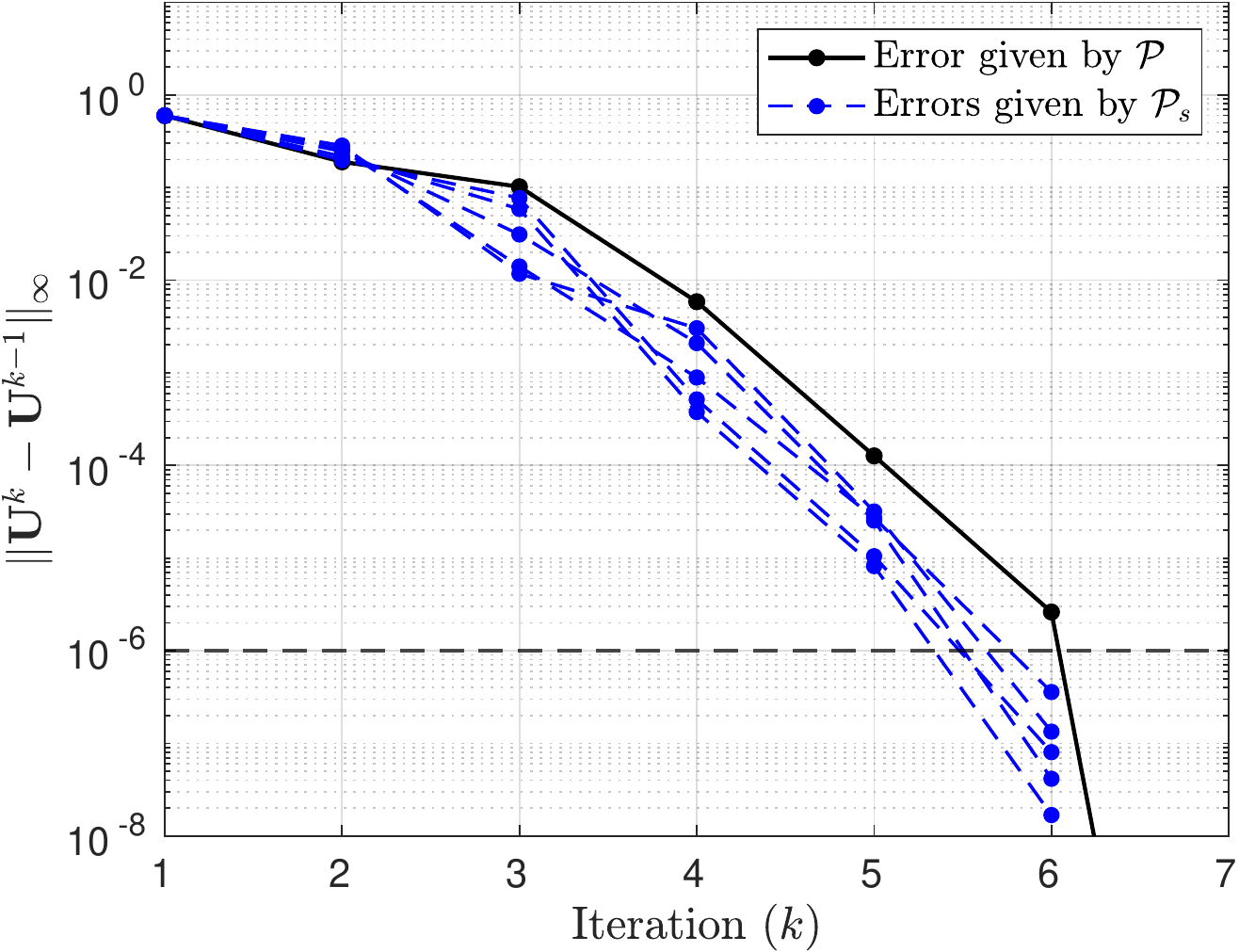} \label{fig:bruss_sol_b}}    
    \caption{(a) Numerical solution of system (4.2) in the phase plane using $\mathcal{F}$ and $\mathcal{P}_s$. Note again that only a subset of the fine times steps of the (mean) $\mathcal{P}_s$ solution are shown for clarity. (b) Errors at successive iterations of $\mathcal{P}$ (black lines) and five independent realisations of $\mathcal{P}_s$ (blue lines). Horizontal dashed black line represents the tolerance $\varepsilon = 10^{-6}$. Both panels run $\mathcal{P}_s$ with sampling rule 1 and $M=10$.}
    \label{fig:bruss_sol}
\end{figure}

\begin{figure}[htbp]
    \centering
    \includegraphics[width=0.49\textwidth]{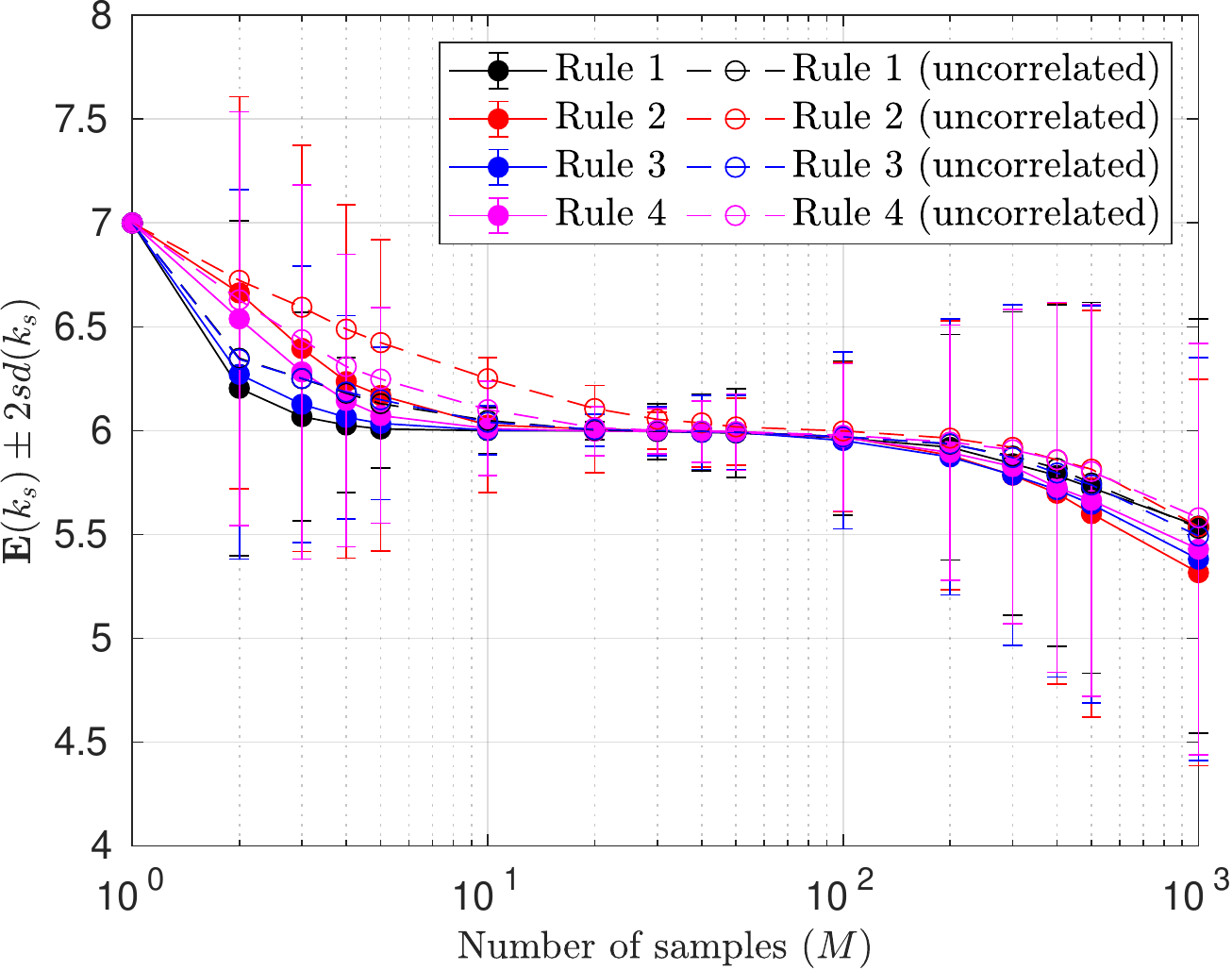}
    \caption{Estimated expectation of $k_s$ as a function of $M$, calculated using estimated distributions of $k_s$ for each sampling rule with error bars representing $\pm$ two standard deviations $sd(k_s)$ shown for correlated (solid lines) sampling rules. Uncorrelated sampling rules are shown with dashed lines without error bars. Distributions are estimated by simulating 2000 independent realisations of $\mathcal{P}_s$ for each $M$.}
    \label{fig:bruss_expec}
\end{figure}

The additional results here complement subsection 4.2. The numerical solutions to system (4.2) in the phase plane are given in \ref{fig:bruss_sol_a} alongside the errors at successive iterations of five runs of $\mathcal{P}_s$ in \ref{fig:bruss_sol_b}. In \cref{fig:bruss_expec} we report the expected value of $k_s$ as a function of $M$ for each of the sampling rules. This result suggests that larger values of $M$ are required to reduce $\mathbb{E}(k_s)$ even further for this stiff system. 

%%%%%%%%%%%%%%%%%%%%%%%%%%%%%%%%%%%%%%%%%%%%%%%%%%%%%%%%%%%%%%%%%
\newpage
\section{`Square limit cycle' system} \label{appen:2Dsquare}
\begin{figure}[t!]
    \centering
    \subfloat[]{\includegraphics[width=0.49\textwidth]{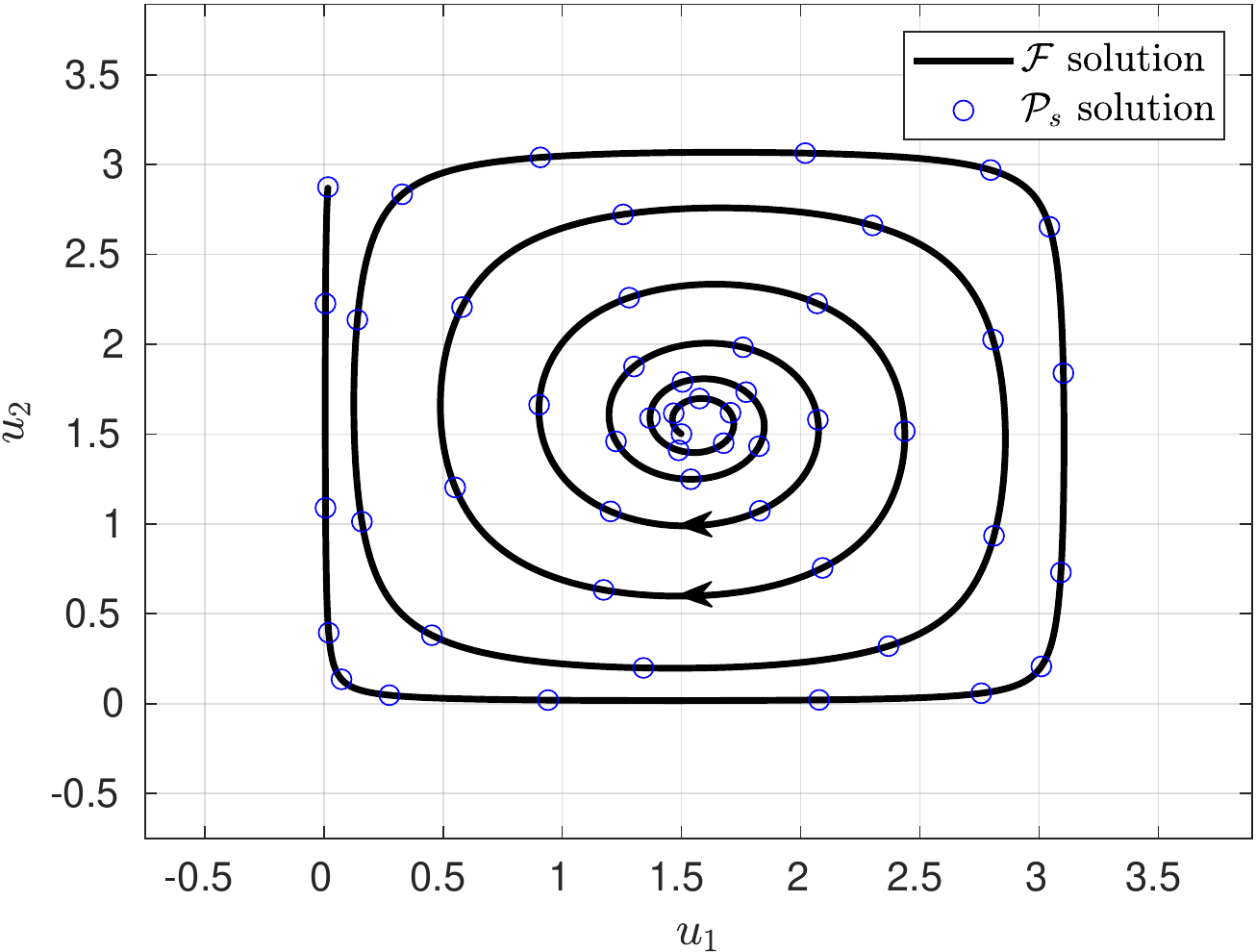} \label{fig:square_sol_a}}
    \subfloat[]{\includegraphics[width=0.49\textwidth]{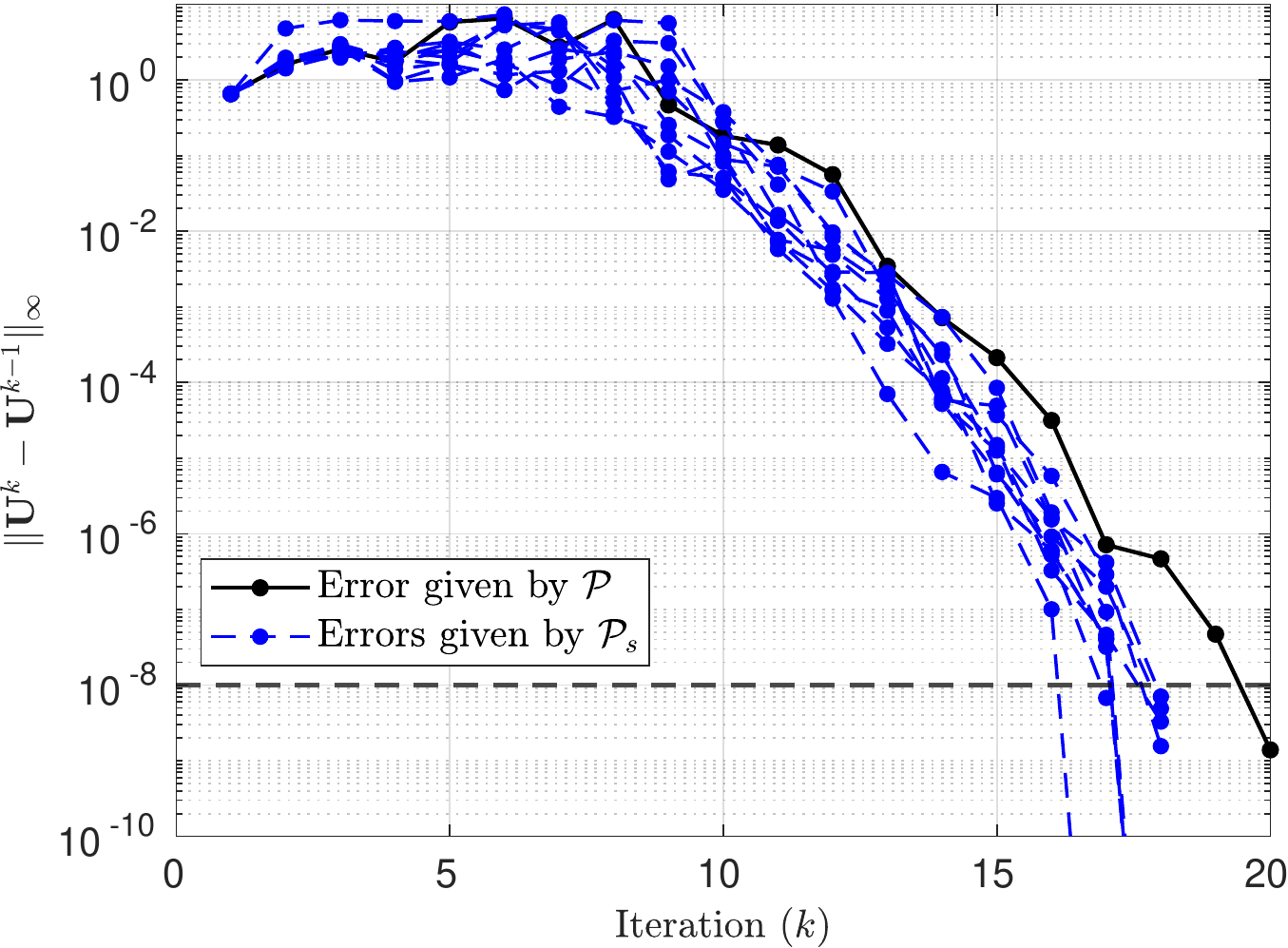} \label{fig:square_sol_b}}      
    \caption{(a) Numerical solution of \eqref{eq:squareODE} over $[0,60]$ using $\mathcal{F}$ serially and $\mathcal{P}_s$. Note again that only a subset of the fine times steps of the (mean) $\mathcal{P}_s$ solution are shown for clarity. (b) Errors at successive iterations of $\mathcal{P}$ and ten independent realisations of $\mathcal{P}_s$. Dashed black line represents the tolerance $\varepsilon = 10^{-8}$. Note that both panels use $\mathcal{P}_s$ with sampling rule 2 and $M=20$.}
    \label{fig:square_sol}
\end{figure}

Consider the system
\begin{subequations} \label{eq:squareODE}
\begin{align} 
        \frac{d u_1}{dt} &= -\sin(u_1) \Big(\frac{\cos(u_1)}{10} + \cos(u_2) \Big), \\
        \frac{d u_2}{dt} &=  -\sin(u_2) \Big(\frac{\cos(u_2)}{10} - \cos(u_1) \Big),
\end{align}    
\end{subequations}
whose solutions, for initial values within the box $[0,\pi] \times [0,\pi]$, converge toward a square-shaped limit cycle on the edges of the box (see \cref{fig:square_sol_a}) \cite{Hirsch2013}. The system is solved on $t \in [0,60]$, starting at $\mathbf{u}(0) = (3/2,3/2)^{\text{T}}$, using $N=30$ sub-intervals and time steps $\delta T = 60/30$ and $\delta t = 60/3000$. As shown in \cref{fig:square_sol_b}, $\mathcal{P}$ takes $k_d = 20$ iterations to converge with tolerance $\varepsilon = 10^{-8}$ whilst the ten realisations of $\mathcal{P}_s$ shown take between $17 \leq k_s \leq 19$.

\begin{figure}[h]
    \centering
    \subfloat[]{\includegraphics[width=0.49\textwidth]{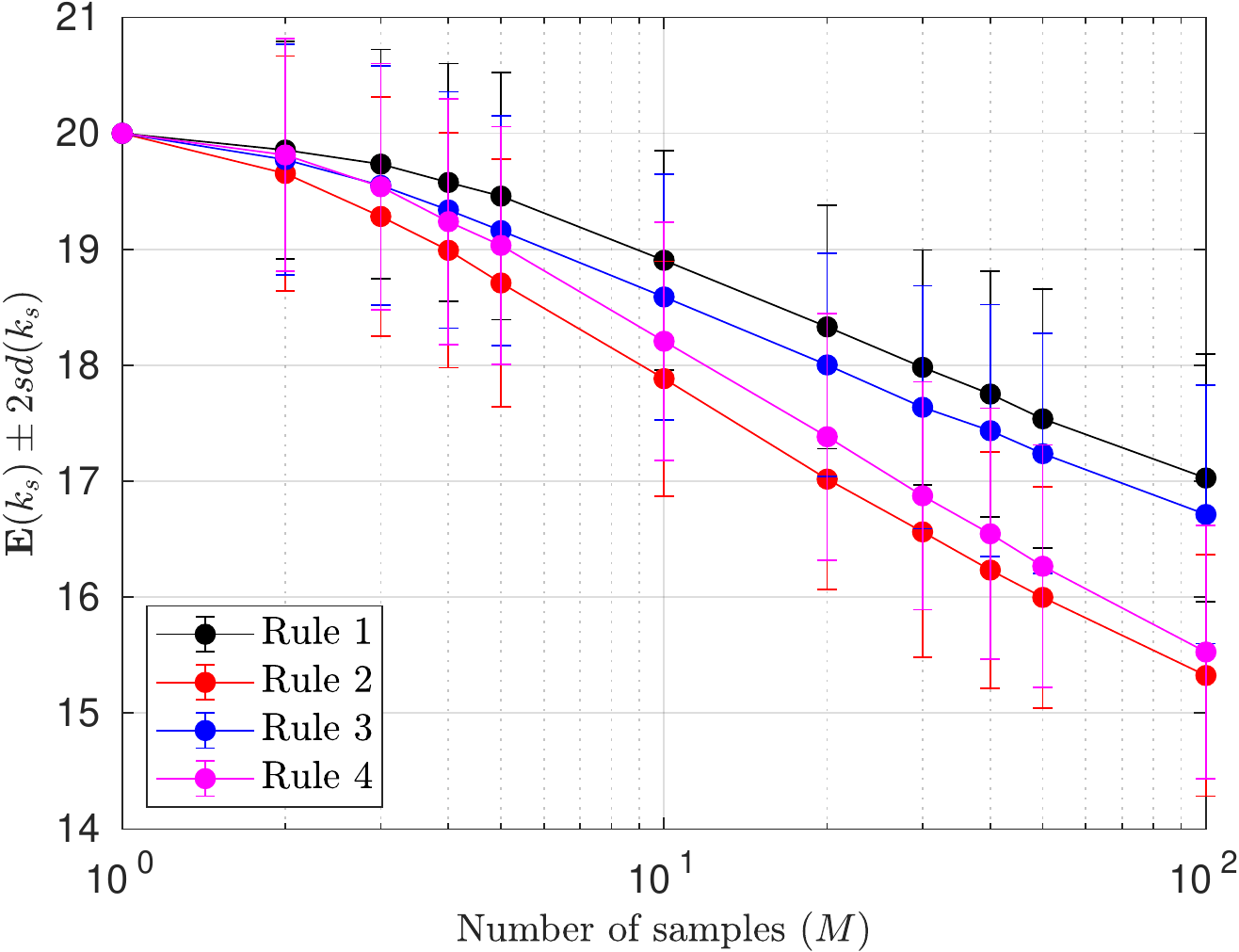} \label{fig:square_probs_a}}
    \subfloat[]{\includegraphics[width=0.49\textwidth]{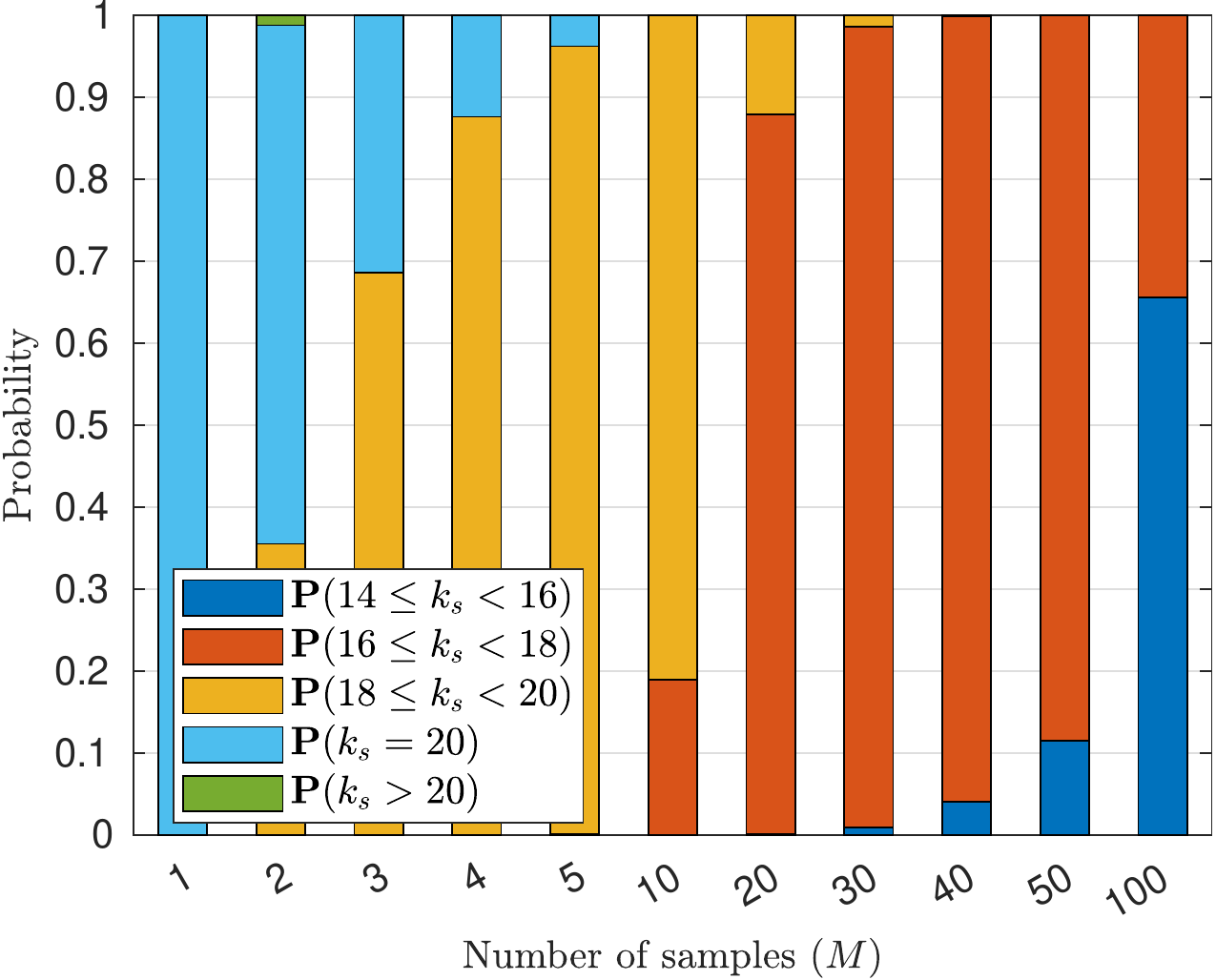} \label{fig:square_probs_b}}
    \caption{(a) Estimated expectation of $k_s$ as a function of $M$, calculated using estimated distributions of $k_s$ for each sampling rule with error bars representing $\pm$ two standard deviations $sd(k_s)$. (b) Estimated discrete distributions of $k_s$ as a function of $M$ for sampling rule 2. Distributions were calculated by simulating 2000 independent realisations of $\mathcal{P}_s$ for each $M$.}
    \label{fig:square_probs}
\end{figure}

Contrary to the previous two-dimensional test problem, \cref{fig:square_probs_a} shows that sampling close to the predictor-corrector values (rules 2 and 4) yield lower expected values of $k_s$. In this case, the bivariate Gaussian outperforms the $t$-copula, however the reverse is true for rules 1 and 3. Results using uncorrelated samples have been shown to generate inferior performance hence are not shown here. The detailed distributions of $k_s$ in \cref{fig:square_probs_b}, using sampling rule 2, show a best performance of $k_s=14$ with 100 samples and $\mathbb{P}(k_s<20)=1$ for $M \approx 10$. They do, however, reveal that in a limited number of cases, using two samples, $\mathcal{P}_s$ can in fact take more than $k_d=20$ iterations to converge ($k_s = 21$ in this small fraction of realisations). This suggests that there may be a minimum number of samples required to beat the convergence rate of parareal for some systems of equations -- something to be investigated in future work.

%%%%%%%%%%%%%%%%%%%%%%%%%%%%%%%%%%%%%%%%%%%%%%%%%%%%%%%%%%%%%%%%%
\newpage
\section{The Lorenz system} \label{appen:3Dlorenz}
\begin{figure}[htbp]
    \centering
    \includegraphics[width=0.49\textwidth]{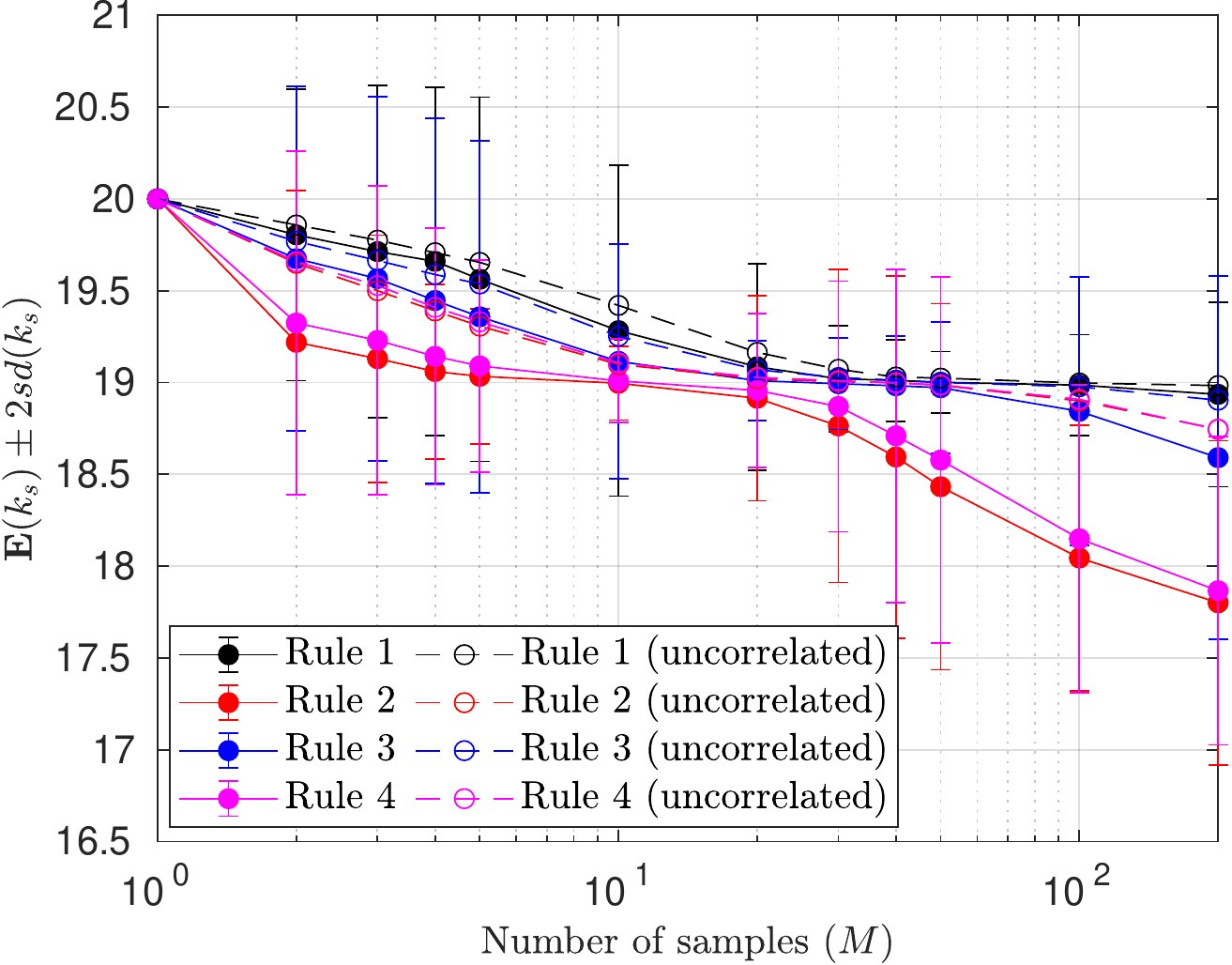}
    \caption{Estimated expectation of $k_s$ as a function of $M$, calculated using estimated distributions of $k_s$ for each sampling rule with error bars representing $\pm$ two standard deviations $sd(k_s)$ shown for correlated (solid lines) sampling rules. Uncorrelated sampling rules are shown with dashed lines without error bars. Distributions are estimated by simulating 2000 independent realisations of $\mathcal{P}_s$ for each $M$.}
    \label{fig:lorenz_expec}
\end{figure}

\begin{figure}[htbp]
    \centering
    \includegraphics[width=0.99\textwidth]{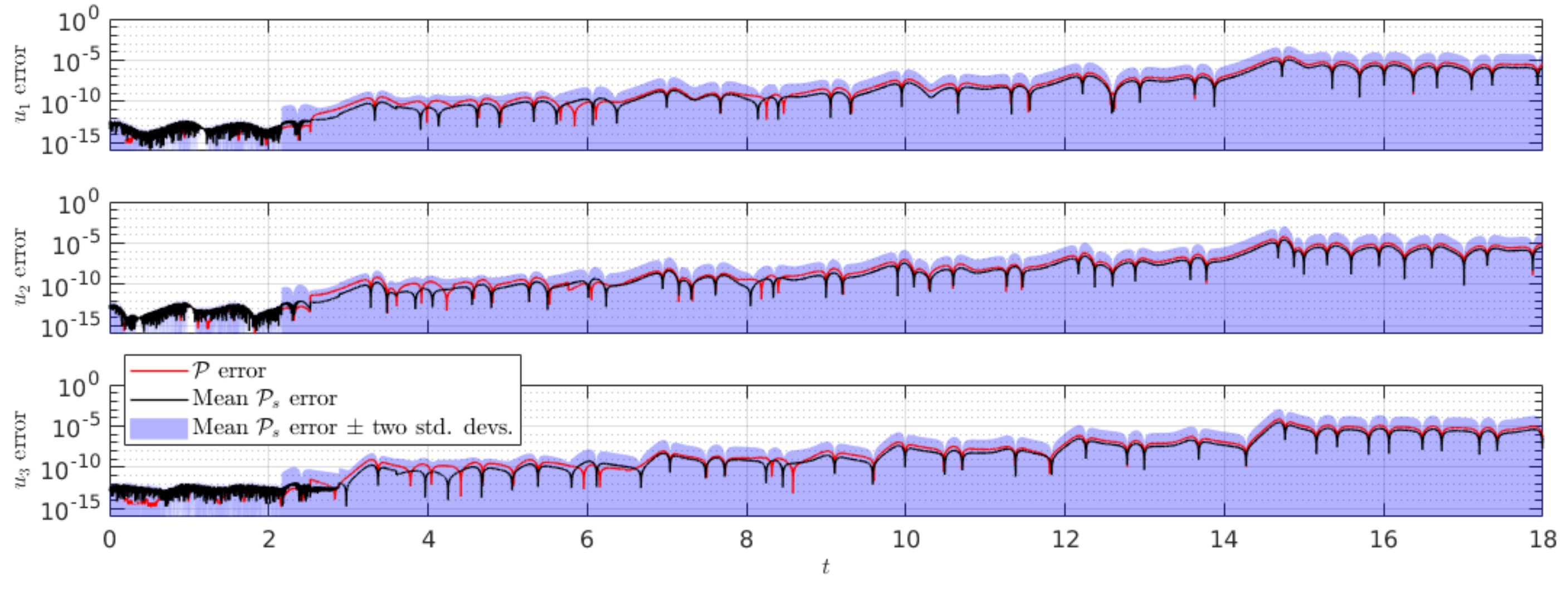}
    \caption{Absolute errors of $\mathcal{P}$ (red) and mean $\mathcal{P}_s$ (black) solutions against the $\mathcal{F}$ solution in each component: $u_1$ (top), $u_2$ (middle), and $u_3$ (bottom). The mean error is obtained by running 2000 independent realisations of $\mathcal{P}_s$ with sampling rule 2 with $M=500$ -- the confidence interval representing the mean $\pm$ two standard deviations is shown in light blue.}
    \label{fig:lorenz_errors}
\end{figure}

The additional results here complement subsection 4.3. \Cref{fig:lorenz_expec} displays $\mathbb{E}(k_s)$ against $M$, indicating that for the Lorenz system, generating correlated samples close to the predictor-corrector solutions (sampling rules 2 and 4) yield the lowest expected values of $k_s$. In \cref{fig:lorenz_errors}, we plot the absolute errors between the mean $\mathcal{P}_s$ solution and the fine solver. Notice that accuracy is maintained with respect to the parareal solution in this chaotic system, even as the errors grow with increasing time.

